\documentclass[12pt,reqno,a4paper]{amsart}
\usepackage{txfonts}

\usepackage{ulem}
\usepackage{mathrsfs}
\usepackage{extarrows}
\usepackage{enumerate}
\usepackage{graphicx}
\usepackage{caption}
\usepackage{float}

\usepackage{color,xcolor}

\usepackage[colorlinks=true,backref=page]{hyperref}
\hypersetup{
    linkcolor=blue,          
    citecolor=red,        
    filecolor=blue,      
    urlcolor=magenta 
}

\usepackage[nobysame,abbrev]{amsrefs}
\renewcommand{\eprint}[1]{{\it Available at}\href{https://arxiv.org/abs/#1}{\it{ arXiv:#1}}.}
\renewcommand{\PrintDOI}[1]{\url{https://doi.org/#1}}
\renewcommand{\MR}[1]{\href{https://mathscinet.ams.org/mathscinet-getitem?mr=#1}{\color{cyan}{MR#1}}}

\BibSpec{article}{%
+{}{\PrintAuthors} {author}
+{,}{ \textrm} {title}
+{:}{ \textrm} {subtitle}
+{.}{ \textit} {journal}
+{,}{ \textbf} {volume}
+{}{ \parenthesize} {date}
+{,}{ } {pages}
+{}{ \parenthesize} {language}
+{.}{} {transition}
+{.}{ \PrintReviews} {review}
+{.}{ \eprint} {eprint}
+{.}{ \PrintDOI} {doi}
}
\BibSpec{book}{%
+{}{\PrintAuthors} {author}
+{,}{ \textit} {title}
+{:}{ \textit} {subtitle}
+{.}{ \textrm} {series} 
+{,}{ Vol.} {volume} 
+{.}{ } {publisher}
+{,}{ } {date}
+{}{ \parenthesize} {language}
+{.}{} {transition}
+{.}{ \PrintReviews} {review}
+{.}{ \PrintDOI} {doi}
}
\usepackage{geometry}
\numberwithin{equation}{section}

\usepackage{thmtools} 
\declaretheorem[name=Theorem, parent=section]{theorem}
\declaretheorem[name=Lemma, sibling=theorem]{lemma}
\declaretheorem[name=Definition, sibling=theorem]{definition}

\declaretheorem[name=Corollary, sibling=theorem]{corollary}
\declaretheorem[name=Remark, sibling=theorem]{remark}
\declaretheorem[name=Example,sibling=theorem]{example}\declaretheorem[name=Proposition,sibling=theorem]{proposition}
\def\bt{\begin{theorem}}
\def\et{\end{theorem}}
\def\bl{\begin{lemma}}
\def\el{\end{lemma}}
\def\bd{\begin{definition}}
\def\ed{\end{definition}}
\def\bp{\begin{proposition}}
\def\ep{\end{proposition}}
\def\bc{\begin{corollary}}
\def\ec{\end{corollary}}
\def\br{\begin{remark}}
\def\er{\end{remark}}
\def\bexa{\begin{example}}
\def\eexa{\end{example}}
\def\bB{{\mathbf B}}
\def\bC{{\mathbf C}}
\def\bL{{\mathbf L}}
\def\bX{{\mathbf X}}
\def\bY{{\mathbf Y}}

\def\1{{\mathbf{1}}}
 
\def\cB{{\mathcal B}}
\def\cM{{\mathcal M}}
\def\cN{{\mathcal N}}
\def\cK{{\mathcal K}}
\def\cP{{\mathcal P}}
\def\cR{{\mathcal R}}

\def\mB{{\mathbb B}}

\def\mD{{\mathbb D}}
\def\mE{{\mathbb E}}

\def\mL{{\mathbb L}}
\def\mN{{\mathbb N}}
\def\mP{{\mathbb P}}

\def\mR{{\mathbb R}}
\def\mS{{\mathbb S}}

\def\sA{{\mathscr A}}
\def\sB{{\mathscr B}}
\def\sE{{\mathscr E}}
\def\sF{{\mathscr F}}
\def\sI{{\mathscr I}}
\def\sJ{{\mathscr J}}
\def\sL{{\mathscr L}}
\def\sM{{\mathscr M}}
\def\sP{{\mathscr P}}
\def\sS{{\mathscr S}}

\def\eps{\varepsilon}

\def\e{\mathrm{e}}

\def\x{{\bf x}}
\def\dif{{\mathord{{\rm d}}}}

\def\p{\partial}
\def\[{{\Big[}}
\def\]{{\Big]}}
\def\<{{\langle}}
\def\>{{\rangle}}
\def\({{\Big(}}
\def\){{\Big)}}

\renewcommand\div{\mathord{{\rm div}\,}}
\def\var{{\mathrm{var}}}

\def\no{\nonumber}
\def\geq{\geqslant}
\def\leq{\leqslant}
\def\ge{\geqslant}
\def\le{\leqslant}
\def\bpf{\begin{proof}}
\def\epf{\end{proof}}

\allowdisplaybreaks
\sloppypar

\begin{document}
	
\title{Supercritical McKean-Vlasov SDE driven by cylindrical $\alpha$-stable process}
\date{\today}
\author{Zimo Hao, Chongyang Ren, and Mingyan Wu}


\address{
Zimo Hao: Fakult\"at f\"ur Mathematik, Universit\"at Bielefeld, 33615, Bielefeld, Germany,
Email: zhao@math.uni-bielefeld.de}

\address{
Chongyan Ren: School of Mathematics and Statistics, Wuhan University, Wuhan, Hubei 430072, P. R. China, Email: rcy.math@whu.edu.cn
}

\address{
Mingyan Wu: School of Mathematical Sciences, Xiamen University, Xiamen, Fujian 361005, P. R. China, Email:  mingyanwu.math@gmail.com }

\thanks{Zimo Hao is supported by the DFG through the CRC 1283 
``Taming uncertainty and profiting from randomness and low regularity in analysis, stochastics and their applications''. Mingyan Wu is partially supported by the National Natural Science Foundation of China (Grant No. 12201227) and the
Fundamental Research Funds for the Central Universities (Grant No. 20720240128). }

\begin{abstract}
In this paper, we study the following supercritical McKean-Vlasov SDE, driven by a symmetric non-degenerate cylindrical $\alpha$-stable process in $\mathbb{R}^d$ with $\alpha \in (0,1)$:
\begin{align*}
\mathord{{\rm d}} X_t = (K * \mu_{t})(X_t)\mathord{{\rm d}}t + \mathord{{\rm d}} L_t^{(\alpha)}, \quad X_0 = x \in \mathbb{R}^d,
\end{align*}
where $K: \mathbb{R}^d \to \mathbb{R}^d$ is a $\beta$-order H\"older continuous function, and $\mu_t$ represents the time marginal distribution of the solution $X$. 
We establish both strong and weak well-posedness under the conditions $\beta \in (1 - \alpha/2, 1)$ and $\beta \in (1 - \alpha, 1)$, respectively. Additionally, we demonstrate strong propagation of chaos for the associated interacting particle system, as well as the convergence of the corresponding Euler approximations. In particular, we prove a commutation property between the particle approximation and the Euler approximation.
\end{abstract}

\maketitle

\section{Introduction}

In recent years, distribution-dependent stochastic differential equations (abbreviated as DDSDEs) have attracted more and more interest, and one important reason is that DDSDEs are deeply related to nonlinear Fokker-Planck equations, which also remain an important area in analysis since they come from many physical and biological models. In the literature, DDSDEs are also often referred to as McKean-Vlasov stochastic differential equations (cf. \cite{Kac56,Kac59,Mc66}). We refer to the article \cite{HRW21} for more about the history and updates on the study of DDSDEs.

\medskip

In this paper, we devote to establishing well-posedness as well as numerical approximations (propagation of chaos and Euler's approximations) for a class of DDSDEs which have H\"older interaction kernels and are driven by symmetric non-degenerate $\alpha$-stable processes with $\alpha \in (0,1)$. Precisely, we consider the following  McKean-Vlasov SDEs in $\mR^d$ with $d\geq 1$:
\begin{align}\label{MV}
\dif X_t=(K*\mu_t)(X_t)\dif t+\dif L_t^{(\alpha)},
\end{align}
where $K:\mR^d\to\mR^d$ is some H\"older continuous function, and $\mu_t$ is the time marginal distribution of the process $(X_t)_{t\ge0}$, and the asterisk denotes the spatial convolution,
\begin{align*}
(K *\mu_t) (x) := \int_{\mR^d} K(x-y) \mu_t(\dif y),
\end{align*}
and $(L_t^{(\alpha)})_{t \geq 0}$ is an $\alpha$-stable process on some probability space $(\Omega,\sF,\mP)$, having a L\'evy measure $\nu^{(\alpha)}$ with the following form:
\begin{align}\label{eq:VV01} 
\nu^{(\alpha)}(A)=\int^\infty_0\left(\int_{\mS^{d-1}}\frac{1_A (r\theta)\Sigma(\dif\theta)}{r^{1+\alpha}}\right)\dif r,\quad A\in\sB(\mR^d),
\end{align} 
where $\Sigma$ is a finite measure over the unit sphere $\mS^{d-1}$ (called spherical measure of $\nu^{(\alpha)}$), and $\sB(\mR^d)$ is the Borel $\sigma$-algebra on $\mR^d$. 

\subsection{ Motivation}\label{Sec:01}

Supposing $\mu_t(\dif x)=\varrho_t(x)\dif x$ in SDE \eqref{MV}, by It\^o's formula, formally, one sees that $\varrho_t$ solves the following nonlocal nonlinear  Fokker-Planck equation in the distributional sense:
\begin{align}\label{in:FPE}
\p_t \varrho_t =(\sL^{(\alpha)})^{*} \varrho_t-\div ((K*\varrho_t) \varrho_t),
\end{align}
where $(\sL^{(\alpha)})^{*}$ is the conjugate operator of the nonlocal operator $\sL^{(\alpha)}$ defined by
\begin{align}\label{eq:nonlocal}
 \sL^{(\alpha)} f (x):=\int_{\mR^d}\( f (x + z) & -f(x)\) \nu^{(\alpha)}(\dif z)
\end{align}
with $\alpha \in (0,1)$ and $f \in C_b^1(\mR^d)$. Notice that when $\nu^{(\alpha)}(\dif z)=c/|z|^{d+\alpha}\dif z$ with an appropriate constant, the operator $\sL^{(\alpha)}$ is exactly the usual fractional Laplacian operator $\Delta^{\alpha/2}:=-(-\Delta)^{\alpha/2}$ on $\mR^d$. 

\medskip

As mentioned at the beginning of this paper, the nonlinear partial differential equation (abbreviated as PDE) \eqref{in:FPE} has attracted extensive and increasing attention due to its relevance to various models from natural science. In a local Laplacian context, for instance, when $K(x)=\nabla \e^{-|x|}$ and $(\sL^{(\alpha)})^{*} \varrho_t$ is substituted with $\div(\varrho_t^2\nabla \varrho_t)$, a so-called aggregation-diffusion equation was established in \cite{TBL06} to model the dynamics of a population where individuals experience long range social attraction and short range dispersal. Considering the case $\sL^{(\alpha)}=\Delta$ (the Laplace operator), when $d=2$ and $K(x)$ is taken as the Biot-Savart law $(-x_2,x_1)/|x|^2$, the equation \eqref{in:FPE} transforms into the famous vorticity form of the 2D Navier-Stokes equation. Moreover, in the same case $\sL^{(\alpha)}=\Delta$, taking $K(x)\asymp |x|^{1-d}$, the equation is seen as the Keller-Segel chemotaxis model  (see \cite{AT21} for example). Back to nonlocal cases, in the fractional Laplace framework, PDE \eqref{in:FPE} reads the surface quasi-geostrophic equation in the mathematical theory of meteorology and oceanography if $K(x)=(-x_2,x_1)/|x|^3$ on $\mR^2$ and $\sL^{(\alpha)}=\Delta^{\alpha/2}$. For further insights into these models and their applications, we refer readers to \cites{BLR11, HRZ23} and the references therein. Moreover, 
in the context of the inviscid case, where the diffusive term $(\sL^{(\alpha)})^{*} \varrho_t$ is omitted, and the interaction kernel $K$ exhibits H\"older continuous at the origin, specifically $K(x)\sim |x|^{\beta}$ with some $\beta\ge0$, the well-posedness of the aggregation equation \eqref{in:FPE} was also studied in the literature, such as \cites{CDFLS} the references therein. See also \cite{BLR11} for the $L^p$ kernels. It should be noted that the initial data $\mu_0$ is always assumed regular in these studies.

\medskip

On the other hand, the McKean-Vlasov SDE \eqref{MV} can be derived from the following large particle systems: 
\begin{align}\label{in:NPS}
\dif X^{N,j}_t=\frac1N\sum_{i=1}^NK(X^{N,j}_t-X^{N,i}_t)\dif t+\dif L^{(\alpha),j}_t,
\end{align}
where $N\in\mN$, and $j=1,2,\cdots,N$, and the interaction kernel $K$ within the systems models attraction and repulsion interactions between each two particles, and $\{L^{(\alpha),j}\}_{j=1}^\infty$ is a family of i.i.d. $\alpha$-stable processes with the same law as $L^{(\alpha)}$, which represents random phenomena. As a matter of fact, back to the 1950s, using $N$-particle systems, Kac (cf. \cite{Kac56}) derived the spactially homogeneous Boltzmann equation, and to this end, the classical notion of propagation of (Kac's) chaos was formalized. We also refer to the Lecture notes \cite{Szn91} and review paper \cite{CD22} for more details of propagation of chaos. It should be noted that in the study of Boltzmann equation, the pure jump process, like $\alpha$-stable processes, play an important role (see \cite{Ta78} for example). Additionally, notice that the scaling factor $\frac{1}{N}$, which is critical in ensuring convergences of these systems as $N\to\infty$, is referred to as the mean-field scaling (for more details, see \cite[Section 1.1]{Ja14}). Consequently, DDSDEs are also termed as mean-field limit SDEs in the literature.

\subsection{Problem statement: supercriticality} 

In this subsection, we discuss the supercriticality of SDEs when $\alpha\in(0,1)$.

\subsubsection{What is supercriticality?}

First of all, let us explain the meaning of supercriticality from the point of view of PDE. Let $\alpha\in(0,2)$. Replacing $(K*\mu_t)(x)$ by a distribution-free drift $b(x)$ in SDE \eqref{MV}, we consider the following classical SDE:
\begin{align}\label{09:0800}
\dif X_t = b(X_t) \dif t + \dif L_t^{(\alpha)},
\end{align}
where the drift term $b$ belongs to the usual H\"older space $\bC^\beta(\mR^d)$ (see subsection \ref{Sec:Be} for the definitions) with $\beta \in (0,1)$. For simplicity, we assume that $\nu^{(\alpha)}(\dif z)=c/|z|^{d+\alpha}\dif z$ in \eqref{eq:VV01}, where $c$ is an appropriate constant such that the infinitesimal generator of $X$ is given by
\begin{align*}
\Delta^{\alpha/2}+b\cdot \nabla. 
\end{align*}

\medskip

On the one hand, it is well-known that the solution $u$ to the evolution equation
\begin{align*}
    \p_t u=\Delta^{\alpha/2} u,~~ t>0,\quad\quad u(0)=\phi,
\end{align*}
can be represented as $u(t,x)=p_t* \phi(x)$, where $p_t(x)$ is the smooth probability density function of the standard $\alpha$-stable process $L^{(\alpha)}$ (see, for example, \cite{WH23}). This implies that, regardless of the singularity in the initial data $\phi$, the solution $u(t)$ becomes smooth for $t>0$. In contrast, for the transport equation
\begin{align}\label{eq:XM0925}
    \p_t u=b\cdot\nabla u,~~t>0,\quad\quad u(0)=\phi,
\end{align}
firstly, there is a well-known counterexample in the dimensional-one case, 
$$
 b(x)= - 2{\rm sign} (x) \sqrt{|x|},
$$
where the drift term $b$ is only H\"older continuous but the Cauchy problem \eqref{eq:XM0925} has infinitely many solutions from any initial condition (cf. \cite{FGP10,FGP18}). Interested readers are also refered to \cite{MS24} for more literature and recent  developments about the non-uniqueness of (stochastic) transport equations. Secondly, considering the simple case, the constant drift term $b\in\mR^d$, one sees that the solution takes the form 
$$
u(t,x)=\phi(x+bt),
$$
which indicates that the singularity of the initial data propagates into the solution when $t>0$. On the other hand, for the fractional Laplace $\Delta^{\alpha/2}$, the following scaling property holds:
\begin{align*}
\Delta^{\alpha/2}f(\lambda x)=\lambda^{\alpha} \Delta^{\alpha/2}f(x),\quad \forall \lambda>0.
\end{align*}

\medskip

Consequently, when $\alpha>1$, $\Delta^{\alpha/2}$ dominates, and the equation is said to be in the subcritical regime; when $\alpha\in(0,1)$, the gradient $b\cdot\nabla$ is of higher order than the fractional Laplacian, leading to the classification of the SDE as supercritical. The critical case corresponds to $\alpha=1$.

\subsubsection{Challenges and methods}

Notably, there is a well-known counterexample presented in \cite{TTW74}, which highlights that it is much more difficult to deal with the supercritical case than with the subcritical case. Rigorously speaking, in the one-dimensional case, Tanaka, Tsuchiya, and Watanabe \cite{TTW74} established the strong well-posedness (strong existence and pathwise uniqueness) to the SDE \eqref{09:0800} for any bounded measurable drift term $b$ in the subcritical regime $\alpha \in (1,2)$, and for any bounded continuous $b$ in the critical case $\alpha =1 $. Furthermore, they demonstrated that in the supercritical case $\alpha \in (0,1)$, the strong well-posedness is ensured when $b$ is a bounded non-decreasing $\bC^\beta$ function with $\beta>1-\alpha$. However, if $\beta+ \alpha <1 $, there exists a $\bC^\beta$ function $b$ for which both pathwise uniqueness and uniqueness in law fail.

\medskip

Samely, for the multi-dimensional case, following the strong well-posedness results obtained by \cite{Pr12} under the subcritical regime conditioned on $\beta>1-\alpha/2$, solving the strong well-posedness problem in the supercritical regime has also proved to be a challenging task. Fortunately, this problem was progressively resolved in \cite{CSZ18} for $\alpha\in(2/3,1)$, and later in \cite{CMP20, HWW23} for $\alpha\in(1/2,1)$. Ultimately, the complete resolution of this problem was achieved by Chen, Zhang, and Zhao in \cite{CZZ21}, in which interested readers can find further details about the well-posedness of SDE \eqref{09:0800}.

\medskip

Moreover, among the results mentioned above, a common approach is applying a priori estimates, such as Schauder's estimates, of the corresponding PDE:
\begin{align*}
    \p_t u=\sL^{(\alpha)}u+b\cdot \nabla u+f.
\end{align*} 
In particular, to handle PDEs in the supercritical regime, the authors of \cite{CZZ21} employed the energy method and Littlewood-Paley theory, where they established a Besov-type a priori estimate. This was subsequently extended to Schauder's estimates in \cite{SX23, Zh21}.

\medskip

In this paper, we consider supercritical $\alpha$-stable DDSDEs with drifts being $\beta$-H\"older continuous and depending on distributions. Especially, the weak and strong well-posedness hold for $\beta>1-\alpha$ and $\beta>1-\alpha/2$ respectively. The key ingredient in this paper is to establish a priori estimates about time regularity of the solution of the PDE (see \autoref{thm:TT} below), which was not studied in the previous works \cite{CZZ21,SX23,Zh21}.

\subsection{Related works and our contributions}

\subsubsection{Well-posedness of DDSDEs}

The first aim of this paper is to establish the well-posedness of the supercritical DDSDE \eqref{MV} (see \autoref{thm:well}). It is well-known that the standard Brownian motion is also seen as a $2$-stable process, and has the infinitesimal generator $\sL^{(\alpha)}$ being the classical Laplace operator $\Delta$ (cf. \cite{Sa99}). Thus, in the following, we first briefly introduce some literature under the Brownian framework, and then discuss the jumped case.

\medskip

To date, in the special case $\alpha=2$, both well-posedness (see \cites{CF22,  MV21,RZ21} for example) and propagation of chaos (e.g. \cite{JW18,La21,HRZ24}) have been extensively investigated. To the best of our knowledge, there are at least four methods to study the well-posedness of DDSDEs, namely, the fixed point argument (see e.g. \cite{Szn91}), the distribution iterations technique (cf. \cite{Wa18}), the approach of Zvonkin's transform and Krylov's estimates (see \cite{RZ21} for example), the bi-coupling method (see e.g. \cite{HRW23}). Since we do not consider the Brownian case in this paper, we only list some literature here, and more details can be found in these works and the references therein. 

\medskip

Compared with the Brownian case, the McKean-Vlasov SDE \eqref{MV} with noise being $\alpha$-stable and with the kernel $K$ being non-Lipschitz has received less attention. In the subcritical case $\alpha\in(1,2)$,  the well-posedness theory of DDSDEs was considered in \cites{BW99,DH23} and in \cite{HRZ23} for second order models. As for the supercritical case $\alpha\in(0,1)$, to the best of our knowledge, there are only two works \cite{FKKM21,DH24} studied well-posedness under the distributional-dependent framework. Regarding the McKean-Vlasov SDEs with H\"older interaction kernels, the authors, Frikha, Konakov, and Menozzi \cite{FKKM21}, established the well-posedness for $\alpha\in(2/3,2)$, which was, recently,  extended to the case $\alpha \in (1/2,1]$ by Deng and Huang \cite{DH24} via a fixed point method. It is worthy pointing out that these two works do not cover the whole supercritical regime $\alpha \in (0,1)$, and only considered the L\'evy measure given by $\nu^{(\alpha)}(\dif z)\asymp\dif z/|z|^{d+\alpha}$, where the L\'evy noise is referred to as a standard $\alpha$-stable process.

\medskip

It is important to note that the standard $\alpha$-stable process differs essentially from the Brownian motion, as the components of the former are not independent. When the components in an $\alpha$-stable process $(L^{(\alpha),1},..,L^{(\alpha),d})$ are independent $1$-dimensional $\alpha$-stable processes, the process will be called the cylindrical $\alpha$-stable process. We point out that cylindrical $\alpha$-stable processes play an important role not only in the $N$-particle system \eqref{in:NPS} but also in the related coupled i.i.d. system as given by 
\begin{align}\label{in:CNP}
\dif X^{j}_t=(K*\mu^j_t)(X^{j}_t)\dif t+\dif L^{(\alpha),j}_t,
\end{align}
which is important in the rest of this paper (see subsubsection \ref{subsubsec:S-chaos}), where the sequence $L^{(\alpha),1},...,L^{(\alpha),N},...$ consists of independent $\alpha$-stable processes. Hence, it is necessary to consider the cylindrical noise for DDSDE \eqref{MV}. On the other hand, obviously,  the L\'evy measure of the cylindrical $\alpha$-stable process has the following form:
$$
\nu^{(\alpha)}(\dif z)\asymp \sum_{i=1}^d\dif z_i/|z_i|^{1+\alpha}\delta_0(\dif \bar z)
$$
with the Dirac measure at zero $\delta_0$ on $\mR^{d-1}$ and $\bar{z}:=(z_1,..,z_{i-1},z_{i+1},...,z_d)$, and is much more singular than the one of the standard $\alpha$-stable process. Moreover, from the point of view of Fourier analysis, the Fourier symbol of the associated operator $\sL^{(\alpha)}=\sum_{j=1}^d \Delta_j^{\alpha/2}$ is given by $\sum_{j=1}^d|\xi_j|^\alpha$, which is notably more singular than $|\xi|^\alpha$, the symbol of the usual fractional Laplacian $\Delta^{\alpha/2}$, since the former one is not differentiable at each axis. 

\medskip
\noindent
{\bf Contribution:} In contrast to \cite{FKKM21,DH24},  we establish the well-posedness of DDSDE \eqref{MV} with a large class of $\alpha$-stable processes including the standard ones and  cylindrical ones (see the condition {\bf (ND)} in the subsection \ref{subsec:main}), under the whole supercritical regime $\alpha\in(0,1)$ (see Theorem \ref{thm:well} below).

\subsubsection{Numerical approximations of DDSDEs}

Following the discussion on well-posedness, it is natural to consider suitable numerical approximations.

\medskip

Recently, Euler-Maruyama approximations for singular SDEs have attracted significant attention (see \cite{HW23} and the references therein for example). In particular, for DDSDEs, there have been several works, such as \cite{WH23} for the subcritical Nemytskii-type case and \cite{HRZ21,Zh19} for the Brownian case, that investigated the following Euler's scheme for McKean-Vlasov SDEs: for $h\in(0,1)$, for $t \in [0,h]$, 
$
X_t^h := X_0 ,
$
and for $t \in [kh, (k + 1)h]$ with $k = 0, 1, 2, \cdots$, 
$$
X_t^h = X_{kh}^h + (K*\mu^h_{kh})  (X^h_{kh}) (t-kh) + L^{(\alpha)}_t - L^{(\alpha)}_{kh}.
$$
Equivalently, the Euler approximation also can be seen as follows: 
\begin{align}\label{in:ES}
X^h_t= X_0 + \int_0^t (K*\mu^h_{\pi_h(s)}) (X^h_{\pi_h(s)}) \dif s+ L^{(\alpha)}_t,
\end{align}
where $\pi_h(t):=[t/h]h$, and $\mu^h_{\pi_h(t)}$ is the distribution of $X^h_{\pi_h(t)}$. Trivially, we have
$$
X_t^h = X_{\pi_h(t)}^h + (K*\mu^h_{\pi_h(t)})(X^h_{\pi_h(t)}) (t-\pi_h(t)) + L^{(\alpha)}_t - L^{(\alpha)}_{\pi_h(t)}.
$$

\medskip

However, it is worth noting that simulating the distribution $\mu^h_t$ is more computationally challenging. Consequently,  considering simulating the empirical measure, we see an Euler's scheme for the $N$-particle systems \eqref{in:NPS}, expressed as
\begin{align}\label{in:ESN}
\dif X^{N,j,h}_t=\frac1N\sum_{i=1}^NK(X^{N,j,h}_{\pi_h(t)}-X^{N,i,h}_{\pi_h(t)})\dif t+\dif L^{(\alpha),j}_t,\ \ j=1,\dots,N.
\end{align}
Fixing the particle number $N$, we see that the process $(X^{N,1,h},\cdots, X^{N,N,h})$ is exactly an Euler approximation of a distributional-free SDE (see \eqref{eq:DTT} for more detials). Now the situation becomes the distribution-free cylindrical case, in which there has been a lot of literature  establishing the weak and strong Euler's convergence  $X^{N,j,h}\to X^{N,j}$ (as $h\to0$), even for the supercritical case (cf.  \cite{BDG22,LZ23}).

\medskip

\begin{figure}[H]
\centering
\includegraphics[scale=1.1]{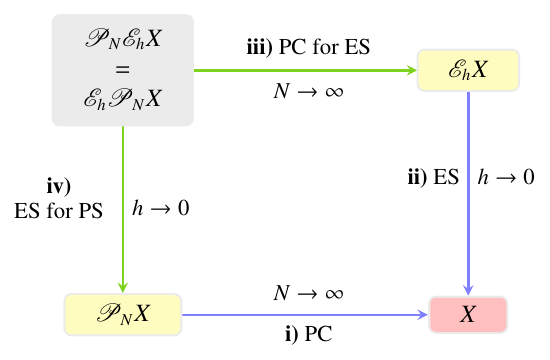}
\caption{Propagation of chaos and Euler's approximations}\label{in:CD}
\end{figure}

Our second aim in this paper is to figure out the following commutative diagram, \autoref{in:CD}, for all $\alpha\in(0,1)$ and $K\in \bC^\beta(\mR^d)$ with $\beta>1-\alpha$.

\noindent
In this context, PC, PS, and ES stand for propagation of chaos, particle system, and Euler's scheme respectively, and we set
$$
\sE_hX:=X^h,\quad \text{and}\quad \sP_NX:=(X^{N,1},...,X^{N,N}),
$$
and
$$
\sP_N \sE_h X= \sE_h \sP_N X:=(X^{N,1,h},\ldots,X^{N,N,h}).
$$
First of all, when the driven noise is a Brownian motion, all of {\bf ii)}, {\bf iii)} and {\bf iv)} in \autoref{in:CD} were studied in \cite{Zh19} for bounded kernel via discrete Krylov's estimates. In addition, the rest part, PC {\bf i)} for singular kernels, also has been investigated well in recent years in both weak and strong sense. Interested readers can find the corresponding results in \cite{JW18,La21,HRZ24}. As for the jumped case, there are fewer works. When the kernel $K$ is Lipschitz and $\alpha\in(1,2)$, the PC {\bf i)} has been obtained in \cite{Gr92}, and the commutative diagram was established in \cite{HY21,FL21}. Moreover, the author in \cite{Ca22} obtained the PC {\bf i)} with $\alpha>1$ and H\"older continuous kernel. The Euler approximation  {\bf iv)} in the strong (path) sense with rate was given by \cite{LZ23} for additive noises under the conditions $K\in \bC^\beta$ with $\beta>1-\alpha/2$ and $\alpha\in(0,2)$ (see  \cite{MX18} for multiplicative case under $\alpha \ge 1$).

\medskip
\noindent
{\bf Contribution:} In this paper, we complete commutative diagram \autoref{in:CD} for the supercritical DDSDEs driven by non-degenerate symmetric $\alpha$-stable processes (e.g. the cylindrical ones, see the condition  {\bf (ND)} in subsection \ref{subsec:main}), where the interaction kernel $K$ exhibits $\beta$-H\"older continuous with $\beta>1-\alpha$. All the convergences are proved both in the weak sense  (see \autoref{thm:appro1} and \autoref{thm:appro1-1} below), and the strong sense (see \autoref{thm:appro2} below). Here, we would like to highlight the significance of the condition $K(x)\sim |x|^{\beta} \in \bC^\beta(\mR^d)$, as it proves valuable in several aggregation models discussed in subsection \ref{Sec:01}.

\subsection{Main results}\label{subsec:main}

Before presenting our main results about supercritical DDSDEs, we introduce some notations and basic assumptions. In the sequel, we denote by $\cP(\mR^d)$ the space consisting of all probability measures on $\mR^d$. Let $L^{(\alpha)}$ be an $\alpha$-stable process whose L\'evy measure $\nu^{(\alpha)}(\dif z)$ given by \eqref{eq:VV01} is symmetric in the sense of $\nu^{(\alpha)}(A) = \nu^{(\alpha)}(-A) $, $\forall A \in \sB(\mR^d)$. Note that for any $\gamma_2>\alpha>\gamma_1\geq 0$, 
\begin{align*} 
\int_{|z|\leq 1} |z|^{\gamma_2} \nu^{(\alpha)}(\dif z) + \int_{|z|>1} |z|^{\gamma_1}\nu^{(\alpha)}(\dif z)<\infty,
\end{align*}
that is
\begin{align}\label{eq:XM02}
\int_{\mR^d} [|z|^{\gamma_1}\wedge |z|^{\gamma_2}]\nu^{(\alpha)}(\dif z)<\infty.
\end{align}
Moreover, throughout this paper, we also assume that the L\'evy measure $\nu^{(\alpha)}$ is symmetric and satisfies the non-degenerate condition:

\medskip
\noindent
$\bf (ND)$ For each $ \theta_0\in \mS^{d-1}$,
$$
\int_{\mS^{d-1}}|\theta\cdot\theta_0|\Sigma(\dif \theta)>0.
$$

\subsubsection{Well-posedness}
Our first goal is to show the well-posedness of the following supercritical DDSDE:
\begin{align}\label{DDSDE}
\dif X_t=b(t,X_t,\mu_t)\dif t+\dif L^{(\alpha)}_t,
\end{align}
where $\mu_t$ is the time marginal distribution of $X_t$, and the drift coefficient $b:\mR_+\times\mR^d\times\cP(\mR^d)\to\mR^d$ satisfies the following assumptions:

\medskip
\noindent
{\bf(H$^\beta_b$)} There is some $\beta\in(0,1)$ and constant $c_0>0$ such that for any $(t,x)\in\mR_+\times\mR^d$ and $\mu, \tilde \mu\in\cP(\mR^d)$,
\begin{align*}
\|b(t,\cdot,\mu)\|_{\bC^\beta}\le c_0,\quad |b(t,x,\mu)-b(t,x,\tilde \mu)|\le c_0\|\mu-\tilde \mu\|_{\beta,\var}.
\end{align*}
Here, the norm $\|\cdot\|_{\bC^\beta}$ is given by \eqref{eq:BM01}, and the norm $\|\cdot\|_{\beta,\var}$, which can be found in \cite{FKKM21}, is defined by \eqref{eq:CZ01}.
\medskip

We first introduce definitions of weak and strong solutions for DDSDE \eqref{DDSDE}.

\begin{definition}[Weak solution]\label{weaksol}
Let $\mu_0$ be a probability measure on $\mR^d$ and $\alpha\in(0, 1)$. We call a stochastic basis $(\Omega, \sF, (\sF_t)_{t\ge0}, \mP)$ together with a pair of $\sF_t$-adapted processes
$(X,L^{(\alpha)})$ defined on it a weak solution of DDSDE \eqref{DDSDE} with initial distribution $\mu_0$, if 
\begin{enumerate}[\rm(i)]
\item $\mu_0=\mP\circ X_0^{-1}$, and $L^{(\alpha)}$ is a $d$-dimensional non-degenerate and symmetric $\alpha$-stable process on the stochastic basis;
\item For each $t>0$, $\mu_t=\mP\circ X_t^{-1}$ and
\begin{align*}
X_t=X_0+\int_0^t b(s,X_s,\mu_s)\dif s+L_t^{(\alpha)}, \, \mP-\text{a.s}.
\end{align*}
\end{enumerate}
\end{definition}

\bd[Strong solution]
Assume that $\alpha \in (0,1)$.  Let $(\Omega, \sF, (\sF_t)_{t\ge0}, \mP)$  be a stochastic basis, and $L^{(\alpha)}$ be a $\sF_t$-adapted $d$-dimensional non-degenerate and symmetric $\alpha$-stable process on the stochastic basis. An $\sF_t$-adapted process $X$ is called a strong solution of DDSDE \eqref{DDSDE}, if for each $t>0$, $\mu_t=\mP\circ X_t^{-1}$ and
\begin{align*}
X_t=X_0+\int_0^t b(s,X_s,\mu_s)\dif s+L_t^{(\alpha)}, \, \mP-\text{a.s}.
\end{align*}
\ed

The proceeding well-posedness theorem is our first main result.

\bt[Well-posedness]\label{thm:well}
Assume that $\alpha \in (0,1)$, $T>0$, and {\bf(H$^\beta_b$)} holds for some $\beta\in(1-\alpha,1)$.  
Then for any $\mu_0\in\cP(\mR^d)$, there is a unique weak solution to DDSDE \eqref{DDSDE}, and a unique martingale solution $\mP\in\sM^{b}_{\mu_0}$ in the sense of \autoref{martsol}. Moreover, if $\beta\in(1-\alpha/2,1)$, the pathwise uniqueness holds, and there is a unique strong solution.
\et

\subsubsection{Numerical approximations}

After establishing the well-posedness of DDSDE, we devote investigating Euler's approximations and propagation of chaos. Fix $\alpha \in (0,1)$ and $T>0$.  In this part, we consider the supercritical McKean-Vlasov SDE \eqref{DDSDE} when the drift coefficient $b$ has the following convolution-type form:
\begin{align*}
b(x,\mu_t) = (K *\mu_t) (x) := \int_{\mR^d} K(x-y) \mu_t(\dif y).
\end{align*} 
Rewriting DDSDE \eqref{DDSDE}, we get the following type DDSDE:
\begin{align}\label{eq:JU02}
\dif X_t=(K*\mu_t)(X_t) \dif t+\dif L_t^{(\alpha)},\ \ t\in (0,T].
\end{align}
In the sequel, we always assume that the kernel $K $ belongs to a H\"older space $\bC^\beta(\mR^d)$ with some $\beta \in (0,1)$, which, obviously, implies that such $b$ satisfies the condition {\bf(H$^\beta_b$)}.
Precisely, we consider the Euler scheme \eqref{in:ES}, and the $N$-particle systems \eqref{in:NPS} with the coupled limit systems \eqref{in:CNP}, and Euler's approximation for the $N$-particle systems \eqref{in:ESN} with the following coupled Euler's scheme:
\begin{align*}
\dif X^{j,h}_t=(K*\mu^{j,h}_{\pi_h(t)})(X^{j,h}_{\pi_h(t)})\dif t+\dif L^{(\alpha),j}_t,
\end{align*}
where $\{L^{(\alpha),j}\}_{j \geq 1}$ is a sequence of independent $\alpha$-stable processes, and $\mu^{j,h}_t$ is the time marginal distribution of $(X^{j,h}_t)_{t \in [0,T]}$. Notice that if $X^{j,h}_0\overset{(d)}{=}X^h_0$, then $\{X^{j,h}_\cdot\}_{j=1}^\infty$ is a family of independent copies of $X^h_\cdot$.

\medskip

Moreover, the particle system \eqref{in:NPS} can be written as an SDE in $\mR^{dN}$: 
\begin{align}\label{eq:DTT}
\dif \bX_t^N = \cK(\bX_t^N)\dif t + \dif \bL_t^{(\alpha),N},
\end{align}
where $\bL_t^{(\alpha),N} := (L_t^{(\alpha),1},\cdots, L_t^{(\alpha),N})$, and $\bX_t^{N} := (X_t^{N,1},\cdots, X_t^{N,N})$, and for $\x =(x^1,\cdots,x^N)\in \mR^{dN}$,
\begin{align}\label{eq:JU93}
\cK^N(\x):=(\cK_1^N(\x),\cdots, \cK_N^N(\x)):= \Bigg( \frac{1}{N}\sum_{i=1}^N K(x^1-x^i) , \cdots,  \frac{1}{N}\sum_{i=1}^N K(x^N-x^i) \Bigg).
\end{align}
Noticing that
$
\| \cK_j^N\|_{\bC^\beta} \leq  2 \|K\|_{\bC^\beta},
$
by \cite[Theorem 1.1]{CZZ21}, we have that there exists a unique weak solution when $\beta >1-\alpha$, and a unique strong solution when $\beta > 1-\alpha/2$.

\medskip

In the following, we divide our results into two situations:
\begin{align*}
 \textbf{the weak sense},\ \  \beta>1-\alpha; 
 \hspace{1em}
\textbf{the strong sense},\ \  \beta>1-\alpha/2.
\end{align*}

\bt[Weak approximations I]\label{thm:appro1}
Let $\alpha\in (0,1)$ and $T>0$. Assume that $K\in\bC^\beta(\mR^d)$ with some $\beta\in(1-\alpha,1)$, and $\mu_0 \in \cP(\mR^d)$.  Then  for any $t \in [0,T]$, the following statements hold:
\begin{enumerate}
\item[\bf i)]  {\rm(Propagation of chaos)}
Suppose that
the law of $(X^{N,1}_0,\cdots,X^{N,N}_0)$ is invariant under any permutation of $\{1,\cdots,N\}$, and for any $1\leq k\leq N$,
\begin{align*}
\mP\circ\(X^{N,1}_0,\cdots,X^{N,k}_0\)^{-1}
\xlongrightarrow{weakly}
\mu_0^{\otimes k},\ \ \text{as}\ \ N\to\infty.
\end{align*}
Then for any $1\leq k\leq N$,
\begin{align*}
\mP\circ\(X^{N,1}_{t},\cdots,X^{N,k}_{t}\)^{-1}
\xlongrightarrow{weakly}
\mu_{t}^{\otimes k}, \ \ \text{as}\ \  N\to\infty,
\end{align*}
where $\mu_t$ is the time marginal distribtuion of $\mu_{[0,T]} \in \sM_{\mu_0}^b$ (see \autoref{martsol}).
\item[\bf ii)] {\rm(Euler's approximation)}
Supposing $X^h_0\overset{(d)}{=}X_0$, we have
$
X^{h}_t \xlongrightarrow{weakly} X_t, \ \ \text{as}\ \ h\to0.
$
\end{enumerate}
\et

\bt[Weak approximations II]\label{thm:appro1-1}
Under the same assumptions of \autoref{thm:appro1},
supposing that $(X^{N,1,h}_0,..,X^{N,N,h}_0)\overset{(d)}{=}(X^{1,h}_0,...,X^{N,h}_0)$, $N\in \mN$, 
consists of i.i.d. random variables with the common law $\mu_0\overset{(d)}{=}X_0$, we have the following statements for any $t\in [0,T]$:
\begin{enumerate}
\item[\bf iii)] {\rm(Propagation of chaos for Euler's scheme)} For any fixed $h\in(0,1)$, 
$$
X^{N,1,h}_t \xlongrightarrow{weakly} X^{1,h}_t, \ \ \text{as}\ \ N\to\infty.
$$
\item[\bf iv)] {\rm(Euler's approximation for $N$-particle systems)} For any fixed $N\in\mN$, 
we have $$
X^{N,1,h}_t\xlongrightarrow{weakly} X^{N,1}_t, \ \ \text{as}\ \ h\to0.
$$
\end{enumerate}
\et   

\begin{proof}[Proof of \autoref{thm:appro1} and \autoref{thm:appro1-1}]
Assertions {\bf i)} and  {\bf iii)} are straightforward from \autoref{thm:wPC} and (ii) of \autoref{thm:PCE}, respectively. The statement {\bf ii)} and {\bf iv)} follow from (i) of \autoref{thm:wEuler} and \eqref{eq:JU01}.
\end{proof} 
   
\bt[Strong approximations]\label{thm:appro2}
Let $\alpha\in (0,1)$, and $K\in\bC^\beta(\mR^d)$ with some $\beta\in(1-\alpha/2,1)$. Assume that for any $N\in \mN$,
$$
(X^{1}_0,\cdots,X^{N}_0) = (X^{N,1}_0,..,X^{N,N}_0)=(X^{N,1,h}_0,..,X^{N,N,h}_0)=(X^{1,h}_0,...,X^{N,h}_0)
$$
consists of i.i.d. random variables with the common law $\mu_0\overset{(d)}{=}X_0$.  Then for any $T>0$ and $p>0$, we have that
\begin{enumerate}[\bf i)]
\item {\rm(Propagation of chaos)} 
$$  \lim_{N\to\infty}\sup_{1\leq j \leq N}\mE\(\sup_{t\in[0,T]}|X^{N,j}_t-X^j_t|^p\)=0;
$$
\item  {\rm(Euler's approximation)} if $X_0=X^h_0$, then
$\lim_{h\to0}\mE\( \sup_{t\in[0,T]}|X^h_t-X_t|^p\)=0$;
\item {\rm(Propagation of chaos for Euler's scheme)} for any fixed $h\in(0,1)$, 
$$
\lim_{N\to\infty}\sup_{1\leq j \leq N}\mE\left(\sup_{t\in[0,T]}|X^{N,j,h}_t-X^{j,h}_t|^p\right)=0;
$$
\item{\rm(Euler's approximation for $N$-particle systems)} for any fixed $N$, 
$$
\lim_{h\to0}\mE\left(\sup_{t\in[0,T]}|X^{N,1,h}_t-X^{N,1}_t|^p\right)=0.
$$
\end{enumerate}
\et   

\begin{proof}
\autoref{thm:SPC}, and (ii) of \autoref{thm:wEuler}, and (i) of \autoref{thm:PCE} below imply {\bf i)}, and {\bf ii)}, and {\bf iii)} respectively. The statement {\bf iv)} is directly from \cite{LZ23} and \eqref{eq:DTT}. 
\end{proof}

\subsection*{Structure of the paper} 

The paper is organized as follows. In Section \ref{Sec:Pre}, we prepare some basic concepts and results of Besov spaces and the norm $\|\cdot\|_{\beta,\var}$. In Section \ref{Sec:03}, we study a class of nonlocal supercritical PDEs and, especailly, obtain some a priori estimates in Besov spaces. In Section \ref{Sec:04}, we establish the well-posedness for the supercritical DDSDE \eqref{DDSDE} by utilizing these a priori estimates and Picard's iteration. Finally, we show proofs of numerical approximations in Section \ref{Sec05}.

\subsection*{Conventions and notations}

Throughout this paper, we use the following conventions and notations: As usual, we use $:=$ as a way of definition. Define $\mN_0:= \mN \cup \{0\}$ and $\mR_+:=[0,\infty)$. The letter $c=c(\cdots)$ denotes an unimportant constant, whose value may change in different places. We use $A \asymp B$ and $A\lesssim B$ (or $A\lesssim_c B$) to denote $c^{-1} B \leq A \leq c B$ and $A \leq cB$, respectively, for some unimportant constant $c \geq 1$. If there is no confusion, we also use $A   \lesssim_\lambda B$ to denote $A \leq c(\lambda) B$ when we want to emphasize that the implicit constant $c$ depends on $\lambda$. 
\begin{itemize}
\item  We denote the space of all bounded and continuous real-valued function on $\mR^d$ by $C_b:=C_b(\mR^d)$. Denote the space of all bounded smooth function with compact support by $C_0^\infty:=C_0^\infty(\mR^d)$. We also define the space consisting of all function $f\in C_b$ with $\lim_{x\to\infty}|f(x)|=0$ by $C_0(\mR^d)$.
 
\item  For every $p\in [1,\infty)$, we denote by $L^p$ the space of all $p$-order integrable functions on $\mR^d$ with the norm denoted by $\|\cdot\|_p$. For $p=\infty$, we set
\begin{align*}
    \|f\|_\infty:=\sup_{x\in\mR^d}|f(x)|.
\end{align*}

\item  For a Banach space $\mB$ and $T>0$, $q\in[1,\infty]$, we denote by
$$
\mL_T^q\mB:= L^q([0,T];\mB),\ \  \mL^q_T:=L^q([0,T]\times \mR^d).
$$
\item We use the convention $\frac{1}{\infty}=0$.
\end{itemize}

\section{Preliminary} \label{Sec:Pre}

\subsection{Besov-H\"older spaces}\label{Sec:Be}

In this subsection, we introduce basic concepts and results of Besov spaces. Let $\sS(\mR^d)$ be the Schwartz space of all rapidly decreasing functions on $\mR^d$, and $\sS'(\mR^d)$ be
the dual space of $\sS(\mR^d)$ called Schwartz generalized function (or tempered distribution) space. Given $f\in\sS(\mR^d)$, 
the Fourier transform $\hat f$ and the inverse Fourier transform  $\check f$ are defined by
$$
\hat f(\xi) :=(2 \pi)^{-d/2}\int_{\mR^d} \e^{-i\xi\cdot x}f(x)\dif x, \quad\xi\in\mR^d,
$$
$$
\check f(x) :=(2 \pi)^{-d/2}\int_{\mR^d} \e^{i\xi\cdot x}f(\xi)\dif\xi, \quad x\in\mR^d.
$$
For every $f\in\sS'(\mR^d)$, the Fourier and the inverse Fourier transform are defined in the following way respectively,
\begin{align*}
\<\hat{f},\varphi\>:=\<f,\hat{\varphi}\>,\ \ \<\check{f},\varphi\>:=\<f,\check{\varphi}\>, \ \  \forall\varphi\in\sS(\mR^d).
\end{align*}
Let $\chi:\mR^{d}\to[0,1]$ be a radial smooth function with
\begin{align*}
\chi(\xi)=
\begin{cases}
1, & \ \  |\xi|\leq 1,\\
0, &\ \ |\xi|>3/2.
\end{cases}
\end{align*}
For $\xi \in \mR^d$, define $\psi(\xi):=\chi(\xi)-\chi(2\xi)$ and
\begin{align*}
\psi_j (\xi) := \begin{cases}
\chi(2\xi), &\ \ \text{for}\ \  j =-1,\\
\psi(2^{-j} \xi),& \ \ \text{for}\ \ j \in \mN_0.
\end{cases}
\end{align*}
Denote $B_r := \{\xi\in \mR^d \mid  |\xi|\leq r\}$ for $r>0$. It is easy to see that $\psi\geq 0$,  supp$\psi\subset B_{3/2}\backslash B_{1/2}$, and
\begin{align}\label{eq:SA00}
\sum_{j=-1}^{k}\psi_j(\xi)=\chi(2^{-k}\xi)\to 1,\ \ \hbox{as}\ \ k\to\infty.
\end{align}
For $j\ge -1$, the block operator $\cR_j$ is defined on $\sS'(\mR^d)$ by
\begin{align*} 
\cR_j f (x):=(\psi_j\hat{f})^{\check\,}(x)=\check\psi_j* f (x).
\end{align*}
By the symmetry of $\psi_j$, we have for any $f\in\sS'(\mR^d)$ and $g \in \sS(\mR^d)$,
\begin{align*}
\int_{\mR^d}\cR_jf(x)g(x)\dif x=\int_{\mR^d}f(x)\cR_jg(x)\dif x, \quad  j\ge-1. 
\end{align*}
The cut-off low frequency operator $S_j$ is defined by 
\begin{align}\label{eq:SA01}
S_j f:= \sum_{k= -1}^{j-1} \cR_k f,
\end{align} 
which, by \eqref{eq:SA00}, derives that
\begin{align}\label{eq:SA01-1}
f = \lim_{j \to \infty}S_j f= \sum_{k\geq -1} \cR_k f,
\end{align} 

Now we state the definition of Besov spaces.
\bd[Besov space]\label{iBesov}
For every $s\in\mR$ and $p,q\in[1,\infty]$, the Besov space $\bB_{p,q}^s(\mR^d)$ is defined by
$$
\bB_{p,q}^s(\mR^d):=\Big\{f\in\sS'(\mR^d)\, \big| \, \|f\|_{\bB^s_{p,q}}:= \[ \sum_{j \geq -1}\left( 2^{s j} \|\cR_j f\|_{p} \right)^q \]^{1/q}  <\infty\Big\}.
$$
If $p=q=\infty$, it is in the sense
$$
\bB_{\infty,\infty}^s(\mR^d):=\Big\{f\in\sS'(\mR^d) \, \big| \, \|f\|_{\bB^s_{\infty,\infty}}:= \sup_{j \geq -1} 2^{s j} \|\cR_j f\|_{\infty} <\infty\Big\}.
$$
\ed 

Recall the following Bernstein's inequality (cf. \cite[Lemma 2.1]{BCD11}).

\bl[Bernstein's inequality]\label{Bern}
For each $k\in\mN_0$, there is a constant $c=c(d,k)>0$ such that for all $j\ge-1$ and $1\leq p_1\leq p_2 \leq \infty$,
\begin{align} \label{S2:Bern}
\|\nabla^k\cR_j f\|_{p_2}  \lesssim_{c}  2^{(k+ d (\frac{1}{p_1}-\frac{1}{p_2}))j}\|\cR_j f\|_{p_1}.
\end{align}
In particular, for any $s \in \mR$ and $1\leq p, q \leq \infty$,
\begin{align*}
\|\nabla^k f\|_{\bB^{s}_{p,q}} \lesssim_{k}  \|f\|_{\bB^{s+k}_{p,q}}.
\end{align*}
\el

We also introduce the following interpolation inequality (cf. \cite[Theorem 2.80]{BCD11}).

\bl[Interpolation inequality]
Let $s_1,s_2\in \mR$ with $s_2>s_1$. For any $p\in [1,\infty]$ and $\theta \in (0,1)$, there is a constant $c=c(s_1,s_2,p,\theta)>0$ such that 
\begin{align}\label{S3:Bern}
\|f\|_{\bB_{p,1}^{\theta s_1+ (1-\theta) s_2} } \lesssim_c \|f\|_{\bB^{s_1}_{p,\infty}}^\theta\|f\|_{\bB^{s_2}_{p,\infty}}^{1-\theta}.
\end{align}
In particular, for any $ s_2> 0 > s_1$,
\begin{align}\label{S9:Bern}
\|f\|_{\infty} \lesssim_c \|f\|_{\bB^{s_1}_{\infty,\infty}}^\theta\|f\|_{\bB^{s_2}_{\infty,\infty}}^{1-\theta},
\end{align}
where $\theta=s_2/(s_2-s_1)$.
\el

Consequently, by \eqref{S3:Bern} and \eqref{S2:Bern}, we have the following embedding relations.

\bc
Assume that $s_1< s_2$, $1\leq p_1\leq p_2 \leq \infty$, and $1\leq r_1\leq r_2 \leq \infty$. Then we have
\begin{align}\label{S4:Bern} 
{\bB_{p_1,\infty}^{s_2}}\hookrightarrow{\bB_{p_1,r_1}^{s_1}}\hookrightarrow  {\bB_{p_2,r_2}^{s_1 - d(\frac{1}{p_1}-\frac{1}{p_2}) }}  
\end{align}
\ec

It is worth discussing here the equivalence between the Besov and H\"older  spaces, which is used in various contexts of this paper without much explanation. For $s>0$, let $\bC^s(\mR^d)$ be the classical $s$-order H\"older space consisting of all measurable functions $f:\mR^d\to\mR$ with
\begin{align}\label{eq:BM01}
\|f\|_{\bC^s}:=\sum_{j=0}^{[s]}\|\nabla^jf\|_\infty+[\nabla^{[s]}f]_{\bC^{s-[s]}}<\infty,
\end{align}
where $[s]$ denotes the greatest integer not more than $s$, and
\begin{align*}
[f]_{\bC^\gamma}:=\sup_{h\in\mR^d}\frac{\|f(\cdot +h ) -f \|_\infty}{|h|^\gamma},~\gamma\in(0,1).
\end{align*}
If $s>0 \text{ and $s\notin\mN$}$, the  following equivalence between $\bB_{\infty,\infty}^s (\mR^d)$ and $ \bC^s  (\mR^d)$ holds: (cf. \cite{Tr92})
\begin{align*}
\|f\|_{\bB_{\infty,\infty}^s}\asymp\|f\|_{\bC^s}.
\end{align*}
However, for any $n\in\mN_0$, we only have one side control that is
$
\|f\|_{\bB^n_{\infty,\infty}}\lesssim\|f\|_{\bC^n}.
$

We also need the following commutator estimate,  which can be found in \cite[Lemma 2.3]{CZZ21}.

\bl[Commutator estimate]\label{lem:WG01}
Let $p,q,r\in[1,\infty]$ with $1/r=1/p+1/q$. For any $\beta \in (0,1)$, there is a constant $c >1$ depending only on $d,p,q,r,\beta$ such that for any $j\geq -1$,
$$
\|[\cR_j, f]g\|_r \leq  c  2^{-\beta j } \| f \|_{\bB_{p,\infty}^\beta} \| g\|_q,
$$
where $[\cR_j,f]g:= \cR_j (fg) - f \cR_j g$. 
\el

\subsection{Metrics on spaces of probability measures}\label{sec:KR}

In this paper, we equip the probability measure space $\cP(\mR^d)$ with the following metric. Fix $\beta\in(0,1)$ and define for $\mu,\tilde \mu\in\cP(\mR^d)$,
\begin{align}\label{eq:CZ01}
\|\mu -\tilde  \mu \|_{\beta,\var}:=\sup_{\|\varphi\|_{\bC^\beta}\le1}\left|\int_{\mR^d}\varphi(y)(\mu (\dif y)-\tilde \mu (\dif y))\right|.
\end{align}

\br
Note that the Banach space $\{f\in\sS'(\mR^d)~;~\|f\|_{\beta,\var}<\infty\}$ with the norm $\|\cdot\|_{\beta,\var}$ is the dual space of the H\"older space $\bC^\beta(\mR^d)=\bB^{\beta}_{\infty,\infty}(\mR^d)$. However, it is not isomorphic with $\bB^{-\beta}_{1,1}$ since $(L^\infty)'\ne L^1$.
\er

\br
The Kantorovich-Rubinstein metric is given by (see \cite[Section 8.3]{Bo07} for more details):
\begin{align*}
\|\mu-\tilde{\mu}\|_{1}:=\sup_{\|\varphi\|_{BL}\le1}\left|\int_{\mR^d}\varphi(y)(\mu-\tilde{\mu}) (\dif y)\right|, \quad\text{where}\quad \|\varphi\|_{BL}:=\|\varphi\|_\infty+\sup_{x\ne y}\frac{|\varphi(x)-\varphi(y)|}{|x-y|}. 
\end{align*}
Notice that the metric $\|\cdot\|_{\beta,\var}$ is stronger than the metric $\|\cdot\|_{1}$ since $\|\varphi\|_{\bC^\beta}\lesssim \|\varphi\|_{BL}$. 
Moreover, we have
\begin{align}\label{eq:JU01}
(\cP(\mR^d),\|\cdot\|_{\var}) \subset  (\cP(\mR^d),\|\cdot\|_{\beta,\var}) \subset (\cP(\mR^d),\|\cdot\|_{1}),
\end{align}
where
\begin{align*}
    \|\mu-\tilde{\mu}\|_{\var}:=\sup_{\|\varphi\|_{\infty}\le1}\left|\int_{\mR^d}\varphi(y)(\mu-\tilde{\mu}) (\dif y)\right|.
\end{align*}
\er

Here we give the following well-known property. 
\bp
Given $\beta\in(0,1)$, the metric $\|\cdot\|_{\beta,\var}$ is equivalent to the weak topology in $\cP(\mR^d)$. More precisely, given any $\{\mu_n\}_{n=1}^\infty\subset \cP(\mR^d)$, then $\lim_{n\to\infty}\|\mu_n-\mu\|_{\beta,\var}=0$ if and only if $\mu_n$ converges to $\mu$ w.r.t. the weak topology. 
\ep
\begin{proof}
For the sufficient part, it follows from \eqref{eq:JU01} that $\lim_{n\to\infty}\|\mu_n-\mu\|_{1}=0$. Then based on \cite[Theorem 8.3.2]{Bo07}, we have $\mu_n$ converges to $\mu$ w.r.t. the weak topology. 

On the contrary, if $\mu_n$ converges to $\mu$ w.r.t. the weak topology, then, by Skorokhod's representation theorem, there exist $\mR^d$-valued random variables $X_n$ and $X$ defined on some probability space $(\Omega,\sF,\mP)$ such that the law of $X_n$ and $X$ are $\mu_n$ and $\mu$ respectively, and $X_n$ converges to $X$ as $n\to\infty$, $\mP$-a.s. Then by the definition, we have
\begin{align*}
    \|\mu_n-\mu\|_{\beta,\var}&=\sup_{\|\varphi\|_{\bC^\beta}\le1}\left|\int_{\mR^d}\varphi(y)(\mu_n-\mu) (\dif y)\right|=\sup_{\|\varphi\|_{\bC^\beta}\le1}|\mE\varphi(X_n)-\mE\varphi(X)|\\
    &\le \mE\left[|X_n-X|^\beta\wedge1\right]\to 0,\quad \text{as $n\to\infty$},
\end{align*}
provided by the dominated convergence theorem. This completes the proof.
\end{proof}

Since $\cP(\mR^d)$ is complete w.r.t. the weak topology, as a result, we have the following complete property, which was also given in \cite[p. 853]{FKKM21}.

\bp\label{pro27}
Given $\beta\in(0,1)$, the space $(\cP(\mR^d),\|\cdot\|_{\beta,\var})$ is complete.
\ep
 
\br
When we consider the finite (signed) measure space $\cM(\mR^d)$ equipped with the norm $\|\cdot\|_{\beta,\var}$, the linear space $(\cM(\mR^d),\|\cdot\|_{\beta,\var})$ is not complete for any $\beta\in(0,1)$. Indeed, if $x_n\to x$ in $\mR^d$, then the Dirac measures $\{\delta_{x_n}\}$ converge to $\delta_x$ in the norm $\|\cdot\|_{\beta,\var}$ since
\begin{align*}
\|\delta_{x_n}-\delta_{x}\|_{\beta,\var}\le |x_n-x|^\beta.
\end{align*}
However, $\|\delta_{x_n}-\delta_x\|_{{\var}}=2$ does not lead to convergence. This observation implies that the norm $\|\cdot\|_{\beta,{\var}}$ is not equivalent to the total variation norm. It is worth noting that the space $(\mathcal{M}(\mathbb{R}^d),\|\cdot\|_{\var})$ constitutes a Banach space with $\|\cdot\|_{\beta,\var}\le\|\cdot\|_{\var}$. Consequently, due to the open mapping theorem, $(\mathcal{M}(\mathbb{R}^d),\|\cdot\|_{\beta,{\var}})$ cannot be a Banach space, and thus it would not be complete.
\er

\subsection{Gronwall's inequality of Volterra-type}
We conclude this section by introducing the following statement of Volterra-type Gronwall's inequality which can be found in \cite[Lemma 2.2]{Zh10}.

\bl\label{lem:g}
Let $f \in L_{loc}^{1}(\mR_+;\mR_+)$ and $T>0$. Assume that for some $\gamma_1,\gamma_2\in[0,1)$ and $c_1,c_2>0$,
$$
f(t) \leq c_1 t^{-\gamma_1} + c_2 \int_0^t (t-s)^{-\gamma_2} f(s) \dif s,\ \ 
t \in (0,T].
$$
Then there is a constant $c_3=c_3(c_2,T,\gamma_1,\gamma_2)>0$ such that for all $t \in (0,T]$,
$$
f(t) \leq c_3 c_1 t^{-\gamma_1}.
$$
\el

\section{Time regularity of supercritical PDEs}\label{Sec:03}

Fix $T>0$. In this section, we consider the following  nonlocal PDE:
\begin{align}\label{eq:PDE}
\p_t u = \sL^{(\alpha)} u-\lambda u +b\cdot\nabla u +f ,~t\in(0,T],\quad u (0)=\varphi\in \bB_{p,\infty}^\beta,
\end{align}
where $\alpha, \beta \in (0,1)$, $\lambda\ge0$, $p \in [1,\infty]$, $ b\in \mL_T^\infty \bB_{p,\infty}^\beta$, and $\sL^{(\alpha)}$ is defined by \eqref{eq:nonlocal}. The proceeding theorem is our main result in this section.

\bt[H\"older regularity]\label{thm:TT}
Let $\alpha,\,\beta \in (0,1)$ and $p\in[2,\infty]$ with $\alpha + \beta >1+d/p$. Assume that $T>0$ and $b \in \mL_T^\infty \bB_{p,\infty}^\beta$. 
Then for each $\lambda\in[0,\infty)$ and $(u(0),f)\in \bB_{p, \infty}^\beta\times \mL_T^\infty\bB_{p, \infty}^\beta$, there exists a unique classical solution $u\in \mL_T^1\bB_{p, 1}^{1}$ to the nonlocal PDE \eqref{eq:PDE} in the sense that for all $(t,x)\in[0,T]\times\mR^d$,
\begin{align*}
u(t,x)=u(0,x)+\int_0^t (\sL^{(\alpha)} -\lambda + b\cdot \nabla) u(s,x)\dif s + \int_0^t f(s,x)\dif s.
\end{align*}
Moreover, for every $\gamma\in[0,\alpha)$, there exist a constant $c>0$ depending on $\gamma, d, p, T, \alpha, \beta,\|b \|_{\mL_T^\infty\bB^\beta_{p,\infty}}$, such that
\begin{align}\label{PDEestimate} 
\|u(t)\|_{\bB_{p, \infty}^{\beta+\gamma}} \leq c \( t^{-\gamma/\alpha} \| u(0) \|_{\bB_{p, \infty}^\beta} +( 1+\lambda)^{-\frac{\alpha-\gamma}{\alpha}}\|f \|_{\mL_T^\infty\bB_{p, \infty}^\beta}\),\ \ t >0.
\end{align}
\et

To prove \autoref{thm:TT}, we need the following four results.
 
\bl[cf. {\cite[Lemma 3.1]{CZZ21}}]\label{lem:KK01}
Assume that $\alpha \in (0,1)$ and $p \in[2,\infty)$. Then for any $f \in \sS'(\mR^d)$ with $\cR_j f \in L^p(\mR^d)$, there is a constant $c  = c(\nu^{(\alpha)}, d,p)>0$ such that for $j=0,1,\cdots,$
\begin{align*}
\int_{\mR^d} |\cR_j f |^{p-2}(\cR_j f ) \sL^{(\alpha)} \cR_j f \dif x \leq -c  2^{\alpha j} \|\cR_j f \|_p^{p},
\end{align*}
and for $j=-1$,
\begin{align*}
\int_{\mR^d} |\cR_{-1}f |^{p-2}(\cR_{-1}f ) \sL^{(\alpha)} \cR_{-1} f \dif x \leq 0.
\end{align*}
\el

\bl[cf. {\cite[Lemma 4.7]{SX23}}]
Assume that $\alpha\in(0,1)$. Let $j\ge0$. Suppose that $\cR_j f\in C_0$ and $|\cR_j f(x_0)|=\|\cR_j f\|_\infty$ with some $x_0\in\mR^d$. Then there is a constant $c>0$ independent of $j$ and $f$ such that 
\begin{align}\label{240911:00}
    \text{\rm sign}(\cR_jf(x_0))\sL^{(\alpha)} \cR_j f(x_0)\le -c2^{\alpha j}\|\cR_j f\|_\infty,
\end{align}
where $\text{\rm sign} (x)$ is the sign function of $x$.
\el

\bl[cf. {\cite[Lemma 3.2]{WZ11}}]
Let $f=f(t,x)$ be a smooth function on $\mR_+\times \mR^d$ with $f(t)\in C_0$ for all $t\ge0$. Then $\p_t\|f(t)\|_\infty$ exists for a.e. $t\ge0$ and
\begin{align}\label{240911:02}
    \p_t \|f(t)\|_\infty=(\p_t f)(t,x_t){\rm sign}(f(t,x_t)),
\end{align}
where $x_t$ is the point such that $x\to f(t,x)$ reaches its norm $\pm \|f(t)\|_\infty$.
\el

By \autoref{lem:WG01}, we have the following commutator estimate.

\bc\label{lem:KK02}
For any $\beta \in (0,1)$, there is a constant $c >1$ such that for any $j\geq -1$ and $ p\in[1,\infty]$,
\begin{align*}
\| [\cR_j , b  ]  \nabla u(t) \|_p \leq c  2^{-\beta j } \| b \|_{ \mL_T^\infty\bB_{p,\infty}^\beta} \| \nabla u(t)\|_\infty.
\end{align*}
\ec

Now we are in a position to give the

\begin{proof}[Proof of \autoref{thm:TT}]
It is well known that the nonlocal PDE \eqref{eq:PDE} has a unique smooth solution if 
$$
b,f\in C_0^\infty(\mR^+\times\mR^d),\, u(0)\in C_0^\infty(\mR^d),
$$
see e.g. \cite{Zh12}. Thus, it suffices to prove the a priori estimate \eqref{PDEestimate}, since the method of modification is applicable (see \cite[Theorem 3.5]{CZZ21}).

\medskip

To obtain the desired estimate, we divide the proof into two steps. 
In the first step, we give the following estimate: there are constants $\kappa_0,\kappa_1,c>0$, which are independent of $j$ and $\lambda$ such that for all $j\ge-1$ and $\lambda\ge0$,
\begin{align}\label{240911:01}
\p_t \| \cR_j u(t)\|_p 
& \leq  -\left( \kappa_0 2^{\alpha j} +\lambda- \kappa_1 \right ) \|\cR_j u(t) \|_p+c 2^{-\beta j } \| \nabla u(t)\|_\infty+\|\cR_j f(t)\|_p.
\end{align}
In the second step, we complete the proof by using the Gronwall's inequality. 

\medskip
\noindent 
{\bf (Step 1)} Applying the block operator $\cR_j$ on both sides of \eqref{eq:PDE}, we have
\begin{align}\label{240911:03}
    \p_t \cR_j u = \sL^{(\alpha)} \cR_j u-\lambda \cR_j u+\cR_j(b\cdot\nabla u)+\cR_j f.
\end{align}
 In the following, we use the energy method (see \cite[Lemma 3.4]{Zh21} or \cite[Lemma 3.1]{ZZ18} for example) to handle the case $p<\infty$, and the method in \cite[Theorem 4.8]{SX23} to handle the case $p=\infty$.
 
\medskip
\noindent
{\bf (Step 1.1) Case $p\in[2,\infty)$}:
 Based on \eqref{240911:03}, for any $p\geq 2$, 
\begin{align*}
\<|\cR_j u|^{p-2}(\cR_j u), \p_t \cR_j u\> 
& = \<|\cR_j u|^{p-2}(\cR_j u) ,\sL^{(\alpha)} \cR_j u \> +\<|\cR_j u|^{p-2}(\cR_j u) ,\cR_j(b\cdot\nabla u)\>\\
&\qquad +\<|\cR_j u|^{p-2}(\cR_j u) ,\cR_j f\>-\lambda\|\cR_j u\|_p^p,
\end{align*}
and then 
\begin{align*}
\frac{\p_t \| \cR_j u\|_p^p}{p}+\lambda\|\cR_j u\|_p^p & = 
\<|\cR_j u|^{p-2}(\cR_j u) ,\sL^{(\alpha)} \cR_j u \> + \< |\cR_j u|^{p-2}(\cR_j u),b\cdot\nabla \cR_ju\>\\
& \qquad+\< |\cR_j u|^{p-2}(\cR_j u), [\cR_j, b ]\nabla u\>
 +\<|\cR_j u|^{p-2}(\cR_j u) ,\cR_j f\> \\
 & := \Lambda^{(1)}_j + \Lambda^{(2)}_j+ \Lambda^{(3)}_j+ \Lambda^{(4)}_j.
\end{align*}
 Now we estimate these four terms in turn.
\begin{itemize}
\item
For $\Lambda^{(1)}_j$, by \autoref{lem:KK01}, one sees that there is a constant $c_1>0$ such that
$$
\Lambda^{(1)}_{-1}(t)\leq 0,\ \ \text{and} \ \ \Lambda^{(1)}_{j} (t) \leq -c_1  2^{\alpha j}\|\cR_j u(t) \|_p^{p}\ \ \text{for}\ \ j \geq 0.
$$
\item For $\Lambda^{(2)}_j$, observing that
\begin{align*}
\Lambda^{(2)}_j  &= \< |\cR_j u|^{p-2}(\cR_j u), (b-S_j b)\cdot\nabla \cR_ju\> +  \< |\cR_j u|^{p-2}(\cR_j u), S_j b \cdot\nabla \cR_ju\>\\
& = \< |\cR_j u|^{p-2}(\cR_j u), (b-S_jb)\cdot\nabla \cR_ju\> - \frac{1}{p} \< S_j  \div b,  |\cR_j u|^{p}
\>,
\end{align*}
where $S_j:=\sum_{k= -1}^{j-1} \cR_k$ (see \eqref{eq:SA01}), and following the proof of \cite[Lemma 3.4]{Zh21}, we obtain that there is a constant $c_2 >0$ such that for all $j \geq -1$,
\begin{align*}
 \Lambda^{(2)}_j (t)\leq c_2  2^{(1-\beta+d/p)j} \|b \|_{\mL_T^\infty\bB^\beta_{p,\infty}}\|\cR_j u(t)\|_p^{p}.
\end{align*}
\item For $\Lambda^{(3)}_j$, by H\"older's inequality and the commutator estimate \autoref{lem:KK02}, we get that there is a constant $c_3 >0$ such that for all $j \geq -1$,
\begin{align*}
\Lambda^{(3)}_j (t)&
\leq \| [\cR_j , b \cdot \nabla ] u(t) \|_p \|\cR_j u\|_p^{p-1}.\\
& \leq c_3  2^{-\beta j } \| b \|_{\mL_T^\infty\bB^\beta_{p,\infty}} \| \nabla u(t)\|_\infty \|\cR_j u(t)\|_p^{p-1}.
\end{align*}
\item For $\Lambda^{(4)}_j$, by H\"older's inequality, we have
$$
\Lambda^{(4)}_j (t)\leq  \|\cR_j f(t)\|_p \|\cR_j u(t)\|_p^{p-1}.
$$
\end{itemize}
Consequently, we obtain that  for all $j \geq -1$,
\begin{align*}
\frac{\p_t \| \cR_j u(t)\|_p^p}{p}+\lambda\|\cR_j u\|_p^p 
& \leq   \left( c_2  2^{(1-\beta+d/p)j}  \|\cR_j u(t)\|_p^p +   c_32^{-\beta j }   \| \nabla u(t)\|_\infty \|\cR_j u(t)\|_p^{p-1} \right ) \|b \|_{\mL_T^\infty\bB^\beta_{p,\infty}}\\
& \qquad - c_1 2^{\alpha j}\|\cR_j u(t) \|_p^{p}\1_{j \geq 0} +   \|\cR_j f(t)\|_p \|\cR_j u(t)\|_p^{p-1},
\end{align*} 
Note that $1-\beta + d/p < \alpha$, and then by Young's inequality we have that
\begin{align*}
  c_2   \|b \|_{\mL_T^\infty\bB^\beta_{p,\infty}}  2^{(1-\beta+d/p)j}\le \frac{c_1}2  2^{\alpha j}+c_4
\end{align*}
with some constant $c_4>0$ independent of $j$.
Then dividing both sides by $\|\cR_j u(t) \|_p^{p-1}$, we get that there are two constants $\kappa_0,\kappa_1>0$ independent of $\lambda$ such that for every $j\ge -1$ and $\lambda \geq 0$,
\begin{align*}
\p_t \| \cR_j u(t)\|_p 
& \leq  -\left( \kappa_0 2^{\alpha j} +\lambda- \kappa_1 \right ) \|\cR_j u(t) \|_p+c_5 2^{-\beta j } \| \nabla u(t)\|_\infty+\|\cR_j f(t)\|_p,
\end{align*}
where the constant $c_5>0$ depends on $\|b \|_{\mL_T^\infty\bB^\beta_{p,\infty}}$. This is exactly the inequality \eqref{240911:01}.

\medskip
\noindent
{\bf (Step 1.2) Case $p=\infty$}: It follows from the Riemann–Lebesgue lemma that $\cR_j u(t)\in C_0$. Then let $x_{t,j}$ be the point so that $\cR_j u(t)$ reaches its maximum at $x_{t,j}$. Without loss of generality, we assume that $\cR_j u(t,x_{t,j})>0$. If not, we replace the function $u$ by $-u$. Thus, by the definition, we have that for any $j \geq -1$,
\begin{align*}
    \sL^{(\alpha)}\cR_j u(t,x_{t,j})=\int_{\mR^d}(\cR_j u(t,x_{t,j}+z)-\cR_j u(t,x_{t,j}))\nu^{(\alpha)}(\dif z)\le 0
\end{align*}
and
\begin{align*}
    b\cdot \nabla \cR_j u(t,x_{t,j})=0,
\end{align*}
which by \eqref{240911:03} and \eqref{240911:00} implies that
\begin{align*}
    \p_t\cR_j u(t,x_{t,j})\le -c2^{-\alpha j}\|\cR_j u(t)\|_\infty\1_{j\ge0}-\lambda \cR_j u(t,x_{t,j})+[\cR_j,b]\nabla u(t,x_{t,j})+\cR_j f(t,x_{t,j}).
\end{align*}
Then in view of \eqref{240911:02}, we have
\begin{align*}
    \p_t\|\cR_j u(t)\|_\infty\le -c2^{-\alpha j}\|\cR_j u(t)\|_\infty\1_{j\ge0}-\lambda \|\cR_j u(t)\|_\infty+\| [\cR_j,b]\nabla u(t)\|_\infty+\|\cR_j f(t)\|_\infty.
\end{align*}
Therefore, by the commutator estimate \autoref{lem:KK02}, we obtain \eqref{240911:01} and finish the proof of {\bf (Step 1)}.

\medskip
\noindent 
{\bf (Step 2)} Let $p \in [2,\infty]$ and use the convention $\frac{1}{\infty} =0$. Based on \eqref{240911:01}, by Gronwall's inequality, one sees that for any $j \geq -1$ and $\lambda\ge 0$,
\begin{align}\label{eq:KK03}
\| \cR_j u(t)\|_p 
\lesssim \e^{-(\kappa_0  2^{\alpha j}+\lambda) t} \|\cR_j u(0) \|_p +  \int_0^t \e^{-(\kappa_0 2^{\alpha j}+\lambda) (t-s)} \( 2^{-\beta j } \| \nabla u(s)\|_\infty+\|\cR_j f(s)\|_p\)  \dif s.
\end{align}
Using the elementary fact that for any real number $n,x >0$,  
$
\e^{-x} \leq n^n  x^{-n},
$
we have that for any $ \gamma\geq 0$,
\begin{align*} 
\e^{-\kappa_0 t 2^{\alpha j} } \lesssim \left[(\kappa_0  t)^{-\gamma/\alpha} 2^{-\gamma j}\right]\wedge 1,
\end{align*}
Thus, back to \eqref{eq:KK03}, we have that
\begin{align*}
\begin{split}
 2^{(\beta + \gamma )j }\| \cR_j u(t)\|_p
&  \lesssim (\kappa_0 t)^{-\gamma/\alpha}\| u(0) \|_{\bB_{p,\infty}^{\beta}}
 +  \int_0^t (\kappa_0(t-s))^{-\gamma/\alpha} \e^{-\lambda (t-s)} \| \nabla u(s)\|_\infty  \dif s \\
& \qquad +2^{\gamma j}\|f\|_{\mL^\infty_T\bB^\beta_{p,\infty}}\int_0^t \e^{-\kappa_0  (t-s)2^{\alpha j}}\e^{-\lambda (t-s)} \dif s,
\end{split}
\end{align*}
where by the change of variables, we have that for any $\gamma\in[0,\alpha)$,
\begin{align*}
\int_0^t \e^{-\kappa_0  (t-s)2^{\alpha j}}\e^{-\lambda (t-s)} \dif s& \le \left(\int_0^t \e^{-\kappa_0  (t-s)2^{\alpha j}}\dif s\right)\wedge \left(\int_0^t\e^{-\lambda (t-s)} \dif s\right)\\
&\lesssim_T 2^{-\alpha j}\wedge (1+\lambda)^{-1}\le (1+\lambda)^{-\frac{\alpha-\gamma}{\alpha}}2^{-\gamma j},
\end{align*}
which implies that 
\begin{align}\label{eq:XZ01} 
 \|u(t)\|_{\bB_{p,\infty}^{\beta+\gamma}} 
&  \lesssim (\kappa_0  t)^{-\gamma/\alpha}   \| u(0) \|_{\bB_{p, \infty}^\beta}+ (1+\lambda)^{-\frac{\alpha-\gamma}{\alpha}}\|f \|_{\mL_T^\infty\bB_{p, \infty}^\beta}  \nonumber\\
&\qquad +  \int_0^t (\kappa_0(t-s))^{-\gamma/\alpha} \e^{-\lambda(t-s)}\| \nabla u(s)\|_\infty   \dif s.
\end{align}

\medskip

 Let $\gamma \in (1-\beta+d/p,\alpha)$. Then combined \eqref{eq:XZ01} with the fact that, for $\alpha>\gamma> 1- \beta+d/p$, 
\begin{align}
\|\nabla u(t)\|_\infty \overset{\eqref{S2:Bern}}{\lesssim} \|u(t)\|_{\bB_{\infty,1}^1}\overset{\eqref{S4:Bern}}{\lesssim} \|  u(t)\|_{\bB_{p,1}^{1+d/p}}\overset{\eqref{S4:Bern}}{\lesssim}  \| u(t)\|_{\bB_{p,\infty}^{\beta+\gamma} }, \label{eq:XZ02}
\end{align}
we have
\begin{align*} 
 \|u(t)\|_{\bB_{p,\infty}^{\beta+\gamma}} 
& \lesssim   (\kappa_0 t)^{-\gamma/\alpha}  \| u(0) \|_{\bB_{p, \infty}^\beta} +(1+\lambda)^{-\frac{\alpha-\gamma}{\alpha}}\|f \|_{\mL_T^\infty\bB_{p, \infty}^\beta}\\
&\quad   +    \kappa_0^{-\gamma/\alpha} \int_0^t (t-s)^{-\frac{\gamma}{\alpha}}  \|u(s) \|_{\bB_{p,\infty}^{\beta+\gamma}} \dif s.
\end{align*}
Hence, for any $\gamma \in (1-\beta+d/p,\alpha)$, we obtain
\begin{align}\label{eq:XZ02-1}
\|u(t)\|_{\bB_{p,\infty}^{\beta+\gamma}} \lesssim  t^{-\gamma/\alpha}\| u(0) \|_{\bB_{p, \infty}^\beta}+(1+\lambda)^{-\frac{\alpha-\gamma}{\alpha}}\|f \|_{\mL_T^\infty\bB_{p, \infty}^\beta},
\end{align}
which is due to the Volterra-type Gronwall's inequality \autoref{lem:g}. 

\medskip

Furthermore, by \eqref{eq:XZ02-1}
 and \eqref{eq:XZ02}, one sees that for each $\tilde\gamma\in(1-\beta+d/p,\alpha)$,
$$
\|\nabla u(s)\|_\infty  \lesssim   s^{-\frac{\tilde\gamma}{\alpha}} \| u(0) \|_{\bB_{p, \infty}^\beta} +\|f \|_{\mL_T^\infty\bB_{p, \infty}^\beta}.
$$
Fix $\tilde \gamma$. Taking this result into account in \eqref{eq:XZ01}, one sees that for any $\gamma\in[0,\alpha)$,
\begin{align*}
 \|u(t)\|_{\bB_{p,\infty}^{\beta+\gamma}} 
&  \lesssim (\kappa_0  t)^{-\gamma/\alpha} 
  \| u(0) \|_{\bB_{p, \infty}^\beta}+(1+\lambda)^{-\frac{\alpha-\gamma}{\alpha}}\|f \|_{\mL_T^\infty\bB_{p, \infty}^\beta} \nonumber\\
&\qquad +  \| u(0) \|_{\bB_{p, \infty}^\beta} \int_0^t (\kappa_0(t-s))^{-\gamma/\alpha}  s^{-\tilde \gamma/\alpha} \dif s\\
& \qquad +\| f \|_{\mL^\infty_T\bB_{p, \infty}^\beta} \int_0^t (\kappa_0(t-s))^{-\gamma/\alpha} \e^{-\lambda(t-s)}\dif s\\
& \lesssim_T t^{-\gamma/\alpha}( 1+ t^{1-\frac{\tilde \gamma}{\alpha}} )  \| u(0) \|_{\bB_{p, \infty}^\beta}+(1+\lambda)^{-\frac{\alpha-\gamma}{\alpha}}\|f \|_{\mL_T^\infty\bB_{p, \infty}^\beta},
\end{align*}
where we used the change of variables in the last inequality. The desired estimate follows because of $t \leq T$. This completes the proof.
\end{proof}

\bc
Assume that $p\in [2,\infty]$, $T>0$, and $\alpha \in (0,1)$. If $\lambda =0$ and $f \equiv  0$ in PDE \eqref{eq:PDE}, then for $d/p +1< \alpha + \beta $ and $\beta\in(0,1)$, there is a constant $c=c(d,T,p,\alpha,\beta,\|b\|_{\mL^\infty \bB_{p,\infty}^\beta})>0$ such that
\begin{align}\label{eq:XY00}
\|\nabla u(t)\|_\infty \lesssim t^{-\frac{1+{d}/{p}-\beta}{\alpha}} \|u(0)\|_{\bB_{p,\infty}^{\beta}}.
\end{align}
\ec

\begin{proof}
We only give the proof of the case $p\neq \infty$, since the case $p=\infty$ is easier and similar. By \eqref{eq:SA01-1} and Berstein's inequality \eqref{S2:Bern},
\begin{align*}
\|\nabla u(t)\|_\infty  \lesssim \|  u(t)\|_{\bB_{\infty,1}^1}\lesssim  \|  u(t)\|_{\bB_{p,1}^{1+d/p}},
\end{align*}
where we used the embedding relation \eqref{S4:Bern} in the last inequality.
Furthermore, for large enough $p$ with $\beta<d/p+1 < \alpha + \beta $, by the interpolation inequality \eqref{S3:Bern}, there are two numbers $\theta \in (0,1)$ and $s_0\in (0,\alpha+\beta)$ such that $(1-\theta) (1+d/p) = s_0$ and
\begin{align*}
 \| u(t)\|_{\bB_{p,1}^{1+d/p}} \lesssim \| u(t)\|_{\bB_{p,\infty}^{0}}^\theta\|u(t)\|_{\bB_{p,\infty}^{s_0}}^{1-\theta} \leq \| u(t)\|_{\bB_{p,\infty}^{s_0}}^\theta\|u(t)\|_{\bB_{p,\infty}^{s_0}}^{1-\theta} = \| u(t)\|_{\bB_{p,\infty}^{s_0}},
\end{align*}
which derives the desired estimate by \autoref{thm:TT}. The proof is finished.
\end{proof}

\section{Well-posedness and Euler's scheme of supercritical DDSDEs}\label{Sec:04}

Fix $\alpha \in (0,1)$. This section is devoted to showing \autoref{thm:well}, the well-posedness of McKean-Vlasov SDE \eqref{DDSDE}. In the sequel, for the $\alpha$-stable process $L^{(\alpha)}$, we denote by $N(\dif r, \dif z)$ the associated Poisson random measure:
$$
N((0,t]\times A)  :=  \sum_{s\in(0,t]} \1_{A}(L_s^{(\alpha)} -  L^{(\alpha)}_{s-}),\ \ A\in \sB(\mR^d\setminus \{0\}) , t >0.
$$
Then, by L\'evy-It\^o's decomposition (cf. \cite[Theorem 19.2]{Sa99}), one sees that
\begin{align*}
L^{(\alpha)}_t =  \lim_{\eps \downarrow 0} \int_0^{t} \int_{\eps<|z|\leq 1} z \widetilde N(\dif r,\dif z) + \int_0^{t} \int_{|z|> 1} z   N(\dif r,\dif z),
\end{align*}
where $ \widetilde{N}(\dif r,\dif z):= N(\dif r,\dif z)-\nu^{(\alpha)}(\dif z)\dif r $ is the compensated Poisson random measure. Under the assumptions $\alpha \in (0,1)$ and $\nu^{(\alpha)}$ being symmetric, we have
$$
L^{(\alpha)}_t = \int_0^{t} \int_{\mR^d} z   N(\dif r,\dif z).
$$

\subsection{Stability of SDEs with distribution-free drifts}
Let $T>0$, and $b_1,b_2\in \mL^\infty_T\bC^\beta$ with some $\beta>1-\alpha$. For each $i=1,2$, it is well-known (see \cite{CZZ21} for example) that there is a unique weak solution $ (\Omega^i, \sF^i,(\sF^i_t)_{t\geq 0}, \mP_i; X^i, L^{(\alpha),i})$ to the following classical SDE:
\begin{align}\label{S4:SDE}
X^i_{t}=x+\int_0^tb_i(r,X^i_{r})\dif r+  L_t^{(\alpha),i},\ \ \forall x \in \mR^d, t \in [0,T].
\end{align}
Furthermore, we have the following stability result.

\bt[Stability estimates]\label{thm:App1}
Let $T>0$ and $\alpha \in (0,1)$. Then for any $\beta\in(1-\alpha,1)$, there is a constant $c>0$ depending on $d,T,\alpha,\beta,\|b_1\|_{\mL^\infty_T\bC^\beta}$ such that for any $t\in[0,T]$,
\begin{align*}
\| \mP_1 \circ(X_t^1)^{-1}-\mP_2 \circ(X_t^2)^{-1}\|_{\beta,\var}\le c \int_0^t (t-r)^{-\frac{1-\beta}{\alpha}} \|b_1(r)-b_2(r)\|_{\infty}\dif r,
\end{align*}
where $\|\cdot\|_{\beta,\var}$ is defined by \eqref{eq:CZ01}.
\et
\begin{proof} 
Fix $t\in(0,T]$. It suffices to estimate 
$
\left|\mE\varphi(X^2_{t})-\mE\varphi(X^1_{t})\right|
$
for any $\varphi\in C^\infty_b(\mR^d)$ with $\|\varphi\|_{\bC^\beta}\leq 1$.  Let $\varphi$ be the terminal condition of the following backward PDE:
\begin{align}\label{eq:GH01}
\p_r u^t+\sL^{(\alpha)} u^t+b_1\cdot\nabla u^t=0,\, r \in [0,t),\quad u^t(t)=\varphi,
\end{align}
where $u^t$ is the shifted function $u^t (r,x): = u(t-r,x)$ and $\sL^{(\alpha)}$ is defined by \eqref{eq:nonlocal}. By It\^{o}'s formula, we have that for $i=1,2$,
\begin{align*}
u^t(t,X^i_t) - u^t(0,X_0^i) =  &
 \int_0^t (\p_r +\sL^{(\alpha)}+b_i  \cdot  \nabla ) u^t (r,X^i_r) \dif r 
 \\
 &+  \int_0^{t} \int_{\mR^d} \( u^t (r,X^i_{r-}+ z) - u^t (r,X^i_{r-}) \) \widetilde N^{i}(\dif r, \dif z).
\end{align*}
Note that the second term of the right hand of the above equality is a martingale. Thus, we obtain  
\begin{align*}
\mE u^t (t,X^i_t ) = &\mE \int_0^t  \(\p_r +\sL^{(\alpha)}  + b_i\cdot  \nabla\)  u^t(r,X^i_r ) \dif r + u^t(0,x)   \\
 \overset{\eqref{eq:GH01}}{=}& \mE \int_0^t  \(( b_i -b_1)\cdot \nabla u^t \) (r,X^i_r ) \dif r+ u^t(0,x)  ,
\end{align*}
which implies that for $i=1$,
\begin{align*}
\mE\varphi(X^1_{t}) =\mE u^t(t,X^1_{t})=u^t(0,x),
\end{align*}
and then, by taking $i=2$,
\begin{align*}
\left |\mE\varphi(X^2_{t})-\mE\varphi(X^1_{t})\right|
& = \left |\mE u^t(t,X_t^2)-\mE u^t (t,X^1_{t})\right|
=\left |\mE\int_0^t
\((b_2-b_1)\cdot\nabla u^t\)(r,X^2_r)\dif r\right|\\
& \leq \sup_{\|\varphi\|_{ \bC^\beta} \leq1} \int_0^t   \|\nabla u^t(r)\|_\infty  \| b_2(r)-b_1(r)\|_\infty\dif r\\
& \overset{\eqref{eq:XY00}}{\lesssim}  \int_0^t (t-r)^{-\frac{1-\beta}{\alpha}} \|b_1(r)-b_2(r)\|_{\infty}\dif r.
\end{align*}
The proof is finished.
\end{proof}

\subsection{Weak and strong solutions}

Now we turn to proving \autoref{thm:well} for the supercritical McKean-Vlasov SDE \eqref{DDSDE}. Note that the condition {\bf(H$^\beta_b$)} provides that there is a constant $c_0$ such that
\begin{align}\label{eq:XY01pre}
\sup_{t \in \mR_+,\mu \in \cP(\mR^d)}\| b(t,\cdot, \mu) \|_{\bC^\beta}  \leq c_0 
\end{align}
and
\begin{align}\label{eq:XY01}
\|b(\cdot,\cdot,\mu^{(1)})-b(\cdot,\cdot,\mu^{(2)})\|_\infty
\le c_0\|\mu^{(1)}-\mu^{(2)}\|_{\beta,\var},
\end{align}
where $\|\cdot\|_{\beta,\var}$ is defined by \eqref{eq:CZ01}.

Now we give the
\begin{proof}[Proof of \autoref{thm:well}] 
({\bf Uniqueness}) Let $X$ and $Y$, whose time marginal distributions are $\mu^X_t$ and $\mu^Y_t$ respectively,  be two solutions of the McKean-Vlasov SDE \eqref{DDSDE} with $X_0\overset{(d)}{=}Y_0$. One sees that if $\mu^X_t=\mu^Y_t$ for all $t\in[0,T]$, then the McKean-Vlasov SDE \eqref{DDSDE} reduces to the classical SDE. Consequently, the weak uniqueness and the strong uniqueness are directly from \cite[Theorem 1.1]{CZZ21} since \eqref{eq:XY01pre} holds.
To this end, in view of   \autoref{thm:App1}, we have that
\begin{align*}
\|\mu^X_t-\mu^Y_t\|_{\beta,\var}&\lesssim_T \int_0^t(t-s)^{-\frac{1-\beta}{\alpha}}\|b(s,\cdot,\mu^X_s)-b(s,\cdot,\mu^Y_s)\|_\infty\dif s\\
&\overset{\eqref{eq:XY01}}{\lesssim} \int_0^t(t-s)^{-\frac{1-\beta}{\alpha}} \|\mu^X_s-\mu^Y_s\|_{\beta,\var} \dif s,
\end{align*}
and then get $\sup_{t\in[0,T]}\|\mu^X_t-\mu^Y_t\|_{\beta,\var}=0$ by Volterra-type Gronwall's inequality, \autoref{lem:g}, since $\alpha>1-\beta>0$.

\medskip\noindent
({\bf Existence}) Next, we show the existence. To this ends, we only need to show the existence of the weak solution for $\beta\in(1-\alpha,1)$. Indeed, once we have the weak existence and uniqueness (the existence and uniqueness of $\mu_t$), when $\beta>1-\alpha/2$, the DDSDE \eqref{DDSDE} becomes to a classical SDE, and the existence of the strong solution is directly from well-known results like \cite{CZZ21}. Set $ X^0_t{=}X_0$ for all $t\in[0,T]$, and let $\{X^n_t\}_{t\in[0,T]}$ be the unique weak solution to the following SDE: (cf. \cite{CZZ21}) 
\begin{align}\label{proof1}
X_t^n=X_0+\int_0^t b(s,X_s^n,\mu^{n-1}_s)\dif s+L_t^{(\alpha)},\, n \geq 1,
\end{align}
where $\mu^{n-1}_t$ is the time marginal distribution of $X^{n-1}_t$.

\medskip\noindent
{\it Step 1.} We first prove that $\{\mu_t^{n}\}_{n \geq 0}$ is a Cauchy sequence in $(\cP(\mR^d),\|\cdot\|_{\beta,\var})$  for every $t \in [0,T]$, that is
\begin{align*}
 \lim_{n\to\infty}\sup_{m\in\mN} \|\mu^{n+m}_t-\mu^{n}_t\|_{\beta,\var} =0.
\end{align*}
Indeed, by  \autoref{thm:App1} and \eqref{eq:XY01}, we obtain that
\begin{align*}
\|\mu^{n+m+1}_t-\mu^{n+1}_t\|_{\beta,\var}&\le\int_0^t(t-s)^{-\frac{1-\beta}{\alpha}}\|b(s,\cdot,\mu^{n+m}_s)-b(s,\cdot,\mu^{n}_s)\|_\infty\dif s\\
& \lesssim\int_0^t(t-s)^{-\frac{1-\beta}{\alpha}} \|\mu^{n+m}_s-\mu^{n}_s\|_{\beta,\var} \dif s,
\end{align*}
which, by H\"older's inequality, derives that for any $p > \frac{\alpha}{\alpha + \beta -1}$,
$$
\|\mu^{n+m+1}_t-\mu^{n+1}_t\|_{\beta,\var}^p \lesssim \int_0^t
 \|\mu^{n+m}_s-\mu^{n}_s\|_{\beta,\var}^p \dif s.
$$
Then by Fatou's lemma, one sees that
\begin{align*}
g(t)\lesssim   \int_0^t g(s)\dif s,
\end{align*}
where
\begin{align*}
g(t):=\limsup_{n\to\infty}\sup_{m\in\mN}  \sup_{s\in[0,t]}\|\mu^{n+m}_s-\mu^{n}_s\|_{\beta,\var}^p.
\end{align*}
Thus, by Volterra-type Gronwall's inequality \autoref{lem:g}, we have $g(t)\equiv0$, which means that $\{\mu_t^{n}\}_{n \geq 0}$ is a Cauchy sequence, and then by  \autoref{pro27}, there is a $\mu_t\in\cP(\mR^d)$ such that
\begin{align*}
\lim_{n\to\infty}\|\mu^{n}_t-\mu_t\|_{\beta,\var}=0, \quad \text{for all $t\in[0,T]$.}
\end{align*}
Moreover, we have
\begin{align}
\lim_{n\to\infty}\sup_{s\in[0,t]}\|\mu^{n}_s-\mu_s\|_{\beta,\var}&=\lim_{n\to\infty}\sup_{s\in[0,t]}\lim_{m\to\infty}\|\mu^{n}_s-\mu^{m}_s\|_{\beta,\var}\no\\
&\le \lim_{n\to\infty}\lim_{m\to\infty}\sup_{s\in[0,t]}\|\mu^{n}_s-\mu^{m}_s\|_{\beta,\var}=0.\label{eq:XY02}
\end{align}

\medskip\noindent
{\it Step 2.} 
It is well-known (cf. \cite{CZZ21}) that, under \eqref{eq:XY01pre}, there is a unique weak solution $X$ for the following SDE:
\begin{align}\label{proof2}
X_t=X_0+\int_0^t b(s,X_s,\mu_s)(X_s)\dif t+L_t^{(\alpha)},
\end{align}
where $(\mu_t)_{t\in[0,T]}$ is the same as the ones in {\it Step 1}. Hence, to establish the existence of weak solutions for SDE \eqref{DDSDE}, it suffices to show that 
$$
\mP\circ(X_t)^{-1}=\mu_t.
$$
 In fact, comparing \eqref{proof1} and \eqref{proof2}, and based on the stability result \autoref{thm:App1}, one sees that
\begin{align*}
\|\mP\circ(X_t)^{-1}-\mP
\circ(X^n_t)^{-1}\|_{\beta,\var}\lesssim \int_0^t(t-s)^{-\frac{1-\beta}{\alpha}}\|\mu^{n-1}_s-\mu_s\|_{\beta,\var}\dif s,
\end{align*}
which deduces that
\begin{align*}
\|\mP\circ(X_t)^{-1}-\mu_t\|_{\beta,\var}&\le\|\mP\circ(X_t)^{-1}-\mP_n \circ(X^n_t)^{-1}\|_{\beta,\var}+ \|\mu_t-\mu^n_t\|_{\beta,\var}\\
&\lesssim \sup_{t \in [0,T]}  \|\mu_t-\mu^n_t\|_{\beta,\var} +  \sup_{t \in [0,T]} \|\mu_t-\mu^{n-1}_t\|_{\beta,\var} \int_0^t(t-s)^{-\frac{1-\beta}{\alpha}}\dif s .
\end{align*}
Taking $n \to \infty$, by \eqref{eq:XY02} and the fact $\alpha + \beta >1$, we get 
$$
\sup_{t \in [0,T]}\|\mP\circ(X_t)^{-1}-\mu_t\|_{\beta,\var} =0.
$$
Hence, we obtain the existence and the proof is complete.
\end{proof}

\subsection{Martingale problems}

Before establishing the martingale problem, we prepare some basic notations. Let $\mD_T := D([0,T];\mR^d)$ be the space of all c$\rm\grave{a}$dl$\rm\grave{a}$g (i.e. right continuous and left limits exist) functions from $[0,T]$ to $\mR^d$. In the following, $\mD_T$ is equipped with Skorokhod topology which makes $\mD_T$ into a Polish space. We use $\omega$ to denote a path in $\mD_T$ 
and  $w_t(\omega)=\omega_t$ to the coordinate process. Let $\cB_t:=\sigma\{w_s, s\leq t\}$ be the natural filtration. Denote all the probability measures over $\mD_T$ by $\cP(\mD_T)$.

Now we proceed to give the definiton and the well-poseness of martingale solutions to DDSDE \eqref{DDSDE}.

\bd[Martingle solutions]\label{martsol}
Let $\mu_0\in\cP(\mR^d)$. A probability measure $\mP\in\cP(\mD_T)$ is called a martingale solution of DDSDE \eqref{DDSDE} with an initial distribution $\mu_0$
if $\mP\circ w_0^{-1}=\mu_0$ and for any $f\in C^2_c(\mR^d)$, the process
\begin{align*}
M^f_t :=f(w_t)-f(w_0)-\int^t_0\Big(\sL^{(\alpha)} +b(s,w_s,\mu_s)\cdot\nabla \Big)f(w_s)\dif s
\end{align*}
is a $\cB_t$-martingale under $\mP$, where $\mu_s:=\mP\circ w^{-1}_s$ and $\sL^{(\alpha)}$ is defined by \eqref{eq:nonlocal}. 
We shall use $\sM^{b}_{\mu_0}$ to denote the set of all martingale 
solutions of DDSDE \eqref{DDSDE} associated with $b$ and the initial distribution ${\mu_0}$.
\ed

Now we are in a position to give 

\begin{proof}[Proof of the well-posedness of martingale solution in \autoref{thm:well}]
(i) Let $(X,L^{(\alpha)})$ be a weak solution of DDSDE \eqref{MV-2} with the initial distribution $\mu_0$ in the sense of \autoref{weaksol}. It follows from It\^o's formula that the distribution of $X$ is a martingale solution of DDSDE \eqref{MV-2} in the sense of \autoref{martsol}. Then the existence is straightforward from the weak well-posedness result in \autoref{thm:well}. 

\medskip\noindent
(ii) Next we show the uniqueness. Let $\mP_1,\mP_2\in \sM^{b}_{\mu_0}$ and denote the time marginal distribution by $\mu^1_t$ and $\mu^2_t$ respectively. Thus, based on \cite[Theorem 2.3]{Kur11}, there are two weak solutions, whose distributions are $\mP_1$ and $\mP_2$, to SDE \eqref{S4:SDE} with $b_i(t,x)=b(t,x,\mu^i_t)$, $i=1,2$. Then by the stability result \autoref{thm:App1}, we have
\begin{align*}
\|\mu^1_t-\mu^2_t\|_{\beta,\var}&\lesssim \int_0^t (t-s)^{-\frac{1-\beta}{\alpha}}\|b_1(s)-b_2(s)\|_\infty\dif s\\
&\lesssim \int_0^t (t-s)^{-\frac{1-\beta}{\alpha}}\|\mu^1_s-\mu^2_s\|_{\beta,\var}\dif s,
\end{align*}
which implies that $\mu^1_t=\mu^2_t$ by  Volterra-type Gronwall's inequality \autoref{lem:g}. Finally, following the standard method (see \cite[p.147, Theorem 6.2.3]{SV06} for example), we have $\mP_1=\mP_2$ and finish the proof. 
\end{proof}

\subsection{Application: Euler's convergence rates}

In this subsection, we consider the following McKean-Vlasov SDE with $b=b(x,\mu_t)$:
\begin{align}\label{eq:DDSDE-XM}
\dif X_t=b(X_t,\mu_t)\dif t+\dif L^{(\alpha)}_t.
\end{align}
Fix $h\in(0,1)$ and consider the following Euler's scheme:
\begin{align}\label{eq:Euler-XM}
\dif X^h_t=b(X^h_{\pi_h(t)},\mu^h_{\pi_h(t)})\dif t+\dif L^{(\alpha)}_t,\quad  X^h_0 \overset{(d)}{=} X_0,
\end{align}
where $\pi_h(t):=[t/h]h$, and $\mu^h_t$ is the time marginal distribution of the process $X^h_t$. In the following, we use the It\^o-Tanaka trick and Zvonkin's transform, respectively, to show the weak and strong Euler's convergence rates. First of all, 
noting that by \cite[Corollary 3.1]{LZ23} for every $\alpha\in(0,2)$ and $q>0$,
\begin{align*}
\mE\left[ |L^{(\alpha)}_{\pi_h(s)}-L^{(\alpha)}_{s}|^q \wedge 1\right]
\lesssim 
\begin{cases}
h^{q/\alpha}, & \text{if}~q\in (0,\alpha),\\
h \log h^{-1},  & \text{if}~q= \alpha,\\
h, & \text{if}~q>\alpha,
\end{cases}
\end{align*}
and $\lim_{ n \to +\infty} \frac{\log n}{n^\epsilon} =0$ (for $\forall \epsilon>0$), one sees that for $\alpha \in (0,1)$,
\begin{align}\label{eq:ZJ01}
\mE\left[ |X^h_{\pi_h(s)}-X^h_{s}|^q \wedge 1\right]
&\le \|b\|_\infty^q h^q+ \mE\left[ |L^{(\alpha)}_{\pi_h(s)}-L^{(\alpha)}_{s}|^q \wedge 1\right] \lesssim h^{q \wedge 1},
\end{align}
which implies that
\begin{align}\label{eq:ZJ02}
\begin{split}
\|\mu^h_{\pi_h(s)}-\mu^h_{s}\|_{\beta,\var}& =\sup_{\|\varphi\|_{\bC^\beta\le1}}|\mE\varphi(X^h_{\pi_h(s)})-\mE\varphi(X^h_{s})|\\
&\le\mE\left[ |X^h_{\pi_h(s)}-X^h_{s}|^\beta\wedge 1\right]\lesssim h^{\beta},
\end{split}
\end{align}
where the implicit constants are independent of the time variable $s$.

The following theorem is our main result in this subsection. 

\bt[Convergence rates for Euler's approximation]\label{thm:wEuler}
Assume that the drift term $b$ satisfies the condition {\bf (H$^\beta_b$)} and $T>0$.
\begin{enumerate}[\rm (i)]
\item (Weak convergence rate)
 If $\beta\in(1-\alpha,1)$, then there is a constant $c>0$ such that for all $h\in(0,1)$,
\begin{align}\label{thm:S601}
\sup_{t\in[0,T]}\|\mu^h_t-\mu_t\|_{\beta,\var}\le c h^\beta .
\end{align}
\item (Strong convergence rate)
If $X_0^h = X_0$ and $\beta\in(1-\alpha/2,1)$, then for any $p\ge 2$, there is a constant $c>0$ such that for all $h\in(0,1)$,
\begin{align*}
\mE\left[ \sup_{t\in[0,T]}|X^h_t-X_t|^p\right] \le c h,.
\end{align*}
\end{enumerate}
\et

\br
Note that when $b(x,\mu_t)=K*\mu_t$ with some $K\in\bC^\beta$,  \autoref{thm:wEuler} implies the statement {\bf ii)} in \autoref{thm:appro1} and \autoref{thm:appro2} directly. Moreover, for fixed $N\in\mN$, applying (i) of \autoref{thm:wEuler} to SDE \eqref{eq:DTT-1}  in $\mR^{Nd}$, the assertion {\bf iv)} of \autoref{thm:appro1-1} follows. Indeed, letting $\sA$ be the projection operator defined by
$\sA(\mathbf{x}):= x^1$ for any $\mathbf{x}=(x^1,\cdots,x^N) \in \mR^{Nd}$, one sees that for any $f \in C_b(\mR^d)$,
\begin{align*}
\lim_{h\to0} \mE f (X_t^{N,h}) = \lim_{h\to0}\mE (f \circ \sA) (\bX_t^{N,h}) \overset{\eqref{thm:S601}}{=} \mE (f \circ \sA) (\bX_t^{N}) = \mE f (X_t^{N}).
\end{align*}
\er
\br
It should be noted that in \cite{BDG22}, the strong convergence was obtained for distribution-free case with $\alpha \ge 2/3$. Therein, the rate can go beyond $1/2$ no matter how small $\beta$ is. A future work may consider whether we can obtain a similar rate for DDSDEs.
\er

Now we give the
\begin{proof}[Proof of \autoref{thm:wEuler}]
Define
\begin{align*}
\tilde b_{\mu^h}^h (s,x):= b(x,\mu_{\pi_h(s)}^h), \ \ \text{and}\ \ \tilde b_{\mu} (s,x):= b(x,\mu_s).
\end{align*}
Note that, by definitions and assumptions,
\begin{align}\label{09:06}
\begin{split}
&| \tilde b^h_{\mu^h} (s,X^h_{\pi_h(s)})-\tilde b_{\mu}(s,X^h_{s})|
\le | b (X^h_{\pi_h(s)},\mu^h_{\pi_h(s)})-   b (X^h_{s},\mu^h_{\pi_h(s)})| \\
&\qquad  +| b (X^h_{s},\mu^h_{\pi_h(s)}) -  b (X^h_{s},\mu^h_{s})|+ | b (X^h_{s},\mu^h_{s}) - b (X^h_{s},\mu_{s})|\\
&\lesssim 
|X^h_{\pi_h(s)}-X^h_{s}|^\beta\wedge 1 +\|\mu^h_{\pi_h(s)}-\mu^h_{s}\|_{\beta,\var}+\|\mu^h_{s}-\mu_{s}\|_{\beta,\var}.
\end{split}
\end{align}
\medskip
\noindent
(i) To use the It\^o-Tanaka trick, we consider the following backward PDE with the terminal term $\varphi\in \bC^\beta(\mR^d)$:
\begin{align*}
\p_s u^t+\sL^{(\alpha)} u^t + \tilde b_{\mu}\cdot \nabla u^t=0,~~s\in(0,t),\quad u^t(t)=\varphi,
\end{align*}
where $u^t$ is the shifted function $u^t (s,x): = u(t-s,x)$ and $\sL^{(\alpha)}$ is defined by \eqref{eq:nonlocal}. By the same argument in the proof of \autoref{thm:App1}, one sees that
\begin{align*}
\mE\varphi(X^h_t)-\mE\varphi(X_t)=\mE\int_0^t \left( \tilde b^h_{\mu^h} (s,X^h_{\pi_h(s)})-\tilde b_{\mu}(s,X^h_{s})\right)\cdot \nabla u^t(s,X^h_{s})\dif s,
\end{align*}
which, combining with \eqref{09:06}, \eqref{eq:ZJ01}, and \eqref{eq:ZJ02}, implies that
\begin{align*}
\|\mu^h_{t}-\mu_{t}\|_{\beta,\var}
& \lesssim \int_0^t \left( h^\beta + \|\mu^h_s -\mu_{s}\|_{\beta,\var} \right) \|\nabla u (t-s)\|_\infty \dif s\\
&\overset{\eqref{eq:XY00}}{\lesssim_T} h^\beta +\int_0^t (t-s)^{-\frac{1-\beta}{\alpha}}\|\mu^h_{s}-\mu_{s}\|_{\beta,\var}\dif s,
\end{align*}
where we used the facts $\alpha+\beta >1$ and $t \leq T$. Hence, by Volterra-type Gronwall's inequality \autoref{lem:g}, we establish the desired estimate.

\medskip
\noindent
(ii) In this part, we prove the Euler's strong convergence rate by Zvonkin's transform and \eqref{thm:S601}.  
Consider the following backward PDE:
\begin{align}\label{eq:XM00}
\p_t u+(\sL^{(\alpha)} -\lambda) u+\tilde b_{\mu}  \cdot (\mathbb{I} + \nabla u) =0,\quad u(T)=0,
\end{align}
where $\lambda>0$. Thanks to $\beta >1-\alpha/2$ and \autoref{thm:TT}, by reversing the time variable, there is an $\eps>0$ small enough and a $\lambda$ large enough such that $\frac\alpha2+\beta>1+\eps$ and
\begin{align}\label{eq:XM01}
\|\nabla u\|_{\mL^\infty_T\bC^{\frac\alpha2+\eps}}\le 1/2.
\end{align}
In the following, we fixed this large number $\lambda$. Defining
\begin{align*}
\Phi(t,x):=x+u(t,x),
\end{align*}
by \eqref{eq:XM01}, one sees that
$$
\frac{1}{2}|x-x'| \leq |\Phi(x)-\Phi(x')|\leq \frac{3}{2}|x-x'|,
$$
which implies that, for each $t>0$, 
\begin{align*}
x\to \Phi(t,x) \text{ is a $C^1$-diffeomorphism on } \mR^d.
\end{align*}
Hence, since \eqref{eq:XM01} and
$$
\nabla \Phi_t^{-1}(x)  = (\mathbb{I}+ \nabla u \circ   \Phi_t^{-1}(x))^{-1} = \sum_{k\geq 0} (-\nabla u \circ \Phi_t^{-1}(x))^k, 
$$
we obtain that
$$
\|\nabla \Phi\|_{\mL^\infty_T}+\|\nabla \Phi^{-1}\|_{\mL^\infty_T}\le 2.
$$
Letting
\begin{align*}
M_t := \int_0^t \int_{\mR^d}
\[ u(s,X_{s-}+z) - u(s,X_{s-}) \] \widetilde{N}(\dif s,\dif z) + \int_0^t \int_{\mR^d}
z {N}(\dif s,\dif z),
\end{align*}
\begin{align*}
M^h_t := \int_0^t \int_{\mR^d}
\[ u(s,X^h_{s-}+z) - u(s,X^h_{s-}) \] \widetilde{N}(\dif s,\dif z) + \int_0^t \int_{\mR^d}
z{N}(\dif s,\dif z),
\end{align*}
and 
\begin{align*}
R^h_t:=\int_0^t \left(  b(X^{h}_{\pi_h(s)},\mu_{\pi_h(s)}^h )-  b(X^{h}_s, \mu_s)\right)\cdot (\mathbb{I} + \nabla u) (s,X^{h}_s)\dif s,
\end{align*}
one sees that, by It\^o's formula and PDE \eqref{eq:XM00}, 
\begin{align*}
\Phi(t,X_t) =  \Phi(0,X_0)  + \lambda\int_0^t u(s,X_s)\dif s+ M_t
\end{align*}
and 
\begin{align*}
\Phi(t,X^{h}_t) =  \Phi(0,X_0^h) + \lambda\int_0^t u(s,X^{h}_s)\dif s + M_t^h + R^h_t.
\end{align*}
Then we deuce that
\begin{align}\label{eq:XM03}
 |X_t^h - X_t|^p&\leq \| {\nabla \Phi^{-1}\|_{\mL_T^\infty}^p} |\Phi(t,X_t^h) - \Phi(t,X_t)|^p\nonumber \\
 &\lesssim \lambda^p  \left|  \int_0^t [u(s,X_s^h) - u(s,X_s)] \dif s\right |^p +  |M_t^h - M_t|^p +  |R^h_t|^p \nonumber 
 \\
 & \lesssim \int_0^t |X_s^h -X_s|^p \dif s  + |M_t^h - M_t|^p +   |R^h_t|^p .
\end{align}
Now we estimate the last two terms in turn. Observing that
\begin{align*}
& \left| \Big(  u(s,X^h_{s-}+z) -  u(s,X^h_{s-}) \Big) - \Big(  u(s,X_{s-}+z) -  u(s,X_{s-}) \Big) \right|\\
&\qquad \le \left |X_{s-}^h - X_{s-}\right | 
\| \nabla u(s,\cdot+z) - \nabla u(s,\cdot) \|_{\mL^\infty}\\
&\qquad \lesssim \left |X_{s-}^h - X_{s-}\right | \( \| \nabla  u \|_{\mL_T^\infty} \wedge \( |z|^{\alpha/2+\eps} \|\nabla u\|_{\mL^\infty_T\bC^{\frac{\alpha}{2}+\eps}}\) \)\\
& \qquad \le  \|\nabla u\|_{\mL^\infty_T\bC^{\frac{\alpha}{2}+\eps}} \left |X_{s-}^h - X_{s-}\right | (1 \wedge  |z|^{\alpha/2+\eps} ),
\end{align*}
and based on Kunita’s inequality (cf. \cite[Theorem 2.11]{Ku04}), one sees that for any $p\geq 2$,
\begin{align}\label{eq:XM05}
\mE\left(\sup_{r\in[0,t]} |M_r^h - M_r|^p\right)  \leq \, & \mE\left[ \( \|\nabla u\|_{\mL^\infty_T\bC^{\frac{\alpha}{2}+\eps}}^2 \int_0^t \int_{\mR^d} \left|X_{s-}-X^{h}_{s-}\right|^2 (1 \land |z|^{\alpha/2+\eps})^2\nu^{(\alpha)}(\dif z) \dif s \)^{p/2} \right]\nonumber\\
&+    \|\nabla u \|_{\mL^\infty_T\bC^{\frac{\alpha}{2}+\eps}}^p \mE \int_0^t \int_{\mR^d} \left|X_{s-}-X^{h}_{s-}\right|^p  (1 \land |z|^{\alpha/2+\eps})^p \nu^{(\alpha)} (\dif z) \dif s \nonumber\\
\lesssim\, & \|\nabla u\|^p_{\mL^\infty_T\bC^{\frac\alpha2+\eps}}\mE\int_0^t \sup_{r\in[0,s]}\left|X_r-X^{h}_r\right|^p\dif s,
\end{align}
where we used the fact \eqref{eq:XM02} in the last inequality. Moreover, by \eqref{09:06}, \eqref{eq:ZJ01}, \eqref{thm:S601}, we have
\begin{align*}
\sup_{r\in[0,t]}|R^h_r| &\lesssim \int_0^t  \left(  |X^h_{\pi_h(s)}-X^h_{s}|^\beta\wedge 1  +h^\beta  \right)(1+ \|\nabla u(s)\|_\infty) \dif s\\
&\overset{\eqref{eq:XM01}}{\lesssim_T} 
\int_0^t  \left(  |X^h_{\pi_h(s)}-X^h_{s}|^\beta\wedge 1  \right)  \dif s  +h^\beta , 
\end{align*}
where we used the facts $\alpha+\beta >1$ and $t \leq T$, and then, by \eqref{eq:ZJ01}, 
\begin{align*}
\mE \left(\sup_{r\in[0,t]}|R^h_r|^p\right) 
\lesssim h^{p\beta}+ T^{p-1} \mE \int_0^t  \left(  |X^h_{\pi_h(s)}-X^h_{s}|^{p\beta}\wedge 1  \right)  \dif s   
\lesssim h^{p\beta}+h^{ (p\beta)\wedge 1} \lesssim h,
\end{align*}
where the last inequality is due to $p\beta>{2(1-\alpha/2)>}1$. Combining all the estimates above, back to \eqref{eq:XM03}, we obtain that
\begin{align*}
 \mE \( \sup_{r\in[0,t]} |X^h_{r}-X_{r}|^p\) \lesssim  \int_0^t  \mE \(\sup_{r \in [0,s]}  |X_r^h - X_r|^p\) \dif s  + h,
\end{align*}
which, by Gronwall's inequality \autoref{lem:g}, implies that
\begin{align*}
\mE\left[\sup_{r\in[0,t]}|X^h_r-X_r|^p\right]
\lesssim h.
\end{align*}
The proof is completed.
\end{proof}

\section{Propagation of chaos for supercritical DDSDEs}\label{Sec05}

Fix $\alpha \in (0,1)$ and $T>0$.  In this section, we consider the supercritical McKean-Vlasov SDE with drift $b$ being convolution-type:
\begin{align}\label{MV-2}
\dif X_t=(K*\mu_t)(X_t) \dif t+\dif L_t^{(\alpha)},\ \ t\in (0,T],
\end{align}
where $\mu_t$ is the time marginal distribution of $X_t$, and $K \in \bC^\beta(\mR^d)$.

\subsection{Propagation of chaos}
For fixed $N \in \mN$, let $(\Omega, \sF, (\sF_s)_{s\ge 0}, \mP, \{L^{(\alpha),j}, X^{N,j} \}_{j=1}^N )$ be a weak solution to the following interacting $N$-particle system for \eqref{MV-2}: 
\begin{align}\label{eq:AC01}
\dif X_t^{N,j} = (K*\eta_t^N)(X_t^{N,j}) \dif t + \dif L^{(\alpha),j}_t,\qquad X_0^{N,j}\overset{(d)}{=} \xi_j^N,
\end{align}
where $\{L^{(\alpha),j}\}_{j=1}^N$ be a sequence of independent $d$-dimensional $\alpha$-stable processes with $\alpha \in (0,1)$, $\{\xi^N_j\}_{1\leq j \leq N}$ are random variables, and $\eta_t^N$ is the empirical distribution measure of $N$-particles $\bX_t^N = (X_t^{N,1},\cdots, X_t^{N,N})$ defined  by
\begin{align}\label{eq:epm}
\eta_t^N(\dif y):= \frac{1}{N} \sum_{i=1}^N \delta_{X_t^{N,i}}(\dif y).
\end{align}
Here, $\delta_x$ stands for the Dirac measure at point $x\in \mR^d$. Observe that by definitions,
\begin{align}\label{eq:DT02}
(K*\eta_t^N)(X_t^{N,j})= \frac{1}{N} \sum_{i=1}^N K(X_t^{N,j}-X_t^{N,i}).
\end{align}

\subsubsection{Weak propagation of chaos}

In this subsection we are going to prove the following weak propagation of chaos for the interacting $N$-particle system \eqref{eq:AC01}.

\bt[Weak propagation of chaos] \label{thm:wPC}
Assume that $K\in\bC^\beta(\mR^d)$ with some $\beta\in(1-\alpha,1)$. For each $N\in\mN$,
let $\xi^N_1,\cdots,\xi^N_N$ be $N$-random variables and $\mu_0\in\cP(\mR^d)$. Suppose that
the law of $(\xi^N_1,\cdots,\xi^N_N)$ is invariant under any permutation of $\{1,\cdots,N\}$, and for any $1\leq k\leq N$,
\begin{align}\label{CC15}
\mP\circ\(\xi^N_1,\cdots,\xi^N_k\)^{-1}\xlongrightarrow{weakly} \mu_0^{\otimes k},\ \ \text{as}\ \ N\to\infty.
\end{align}
Then for any $1\leq k\leq N$ and $T>0$,
\begin{align}\label{CC14}
\mP\circ\(X^{N,1}_{[0,T]},\cdots,X^{N,k}_{[0,T]}\)^{-1}\xlongrightarrow{weakly} \mu_{[0,T]}^{\otimes k}, \ \ \text{as}\ \  N\to\infty,
\end{align}
where $\mu_{[0,T]}$ is the law of the unique martingale solution of DDSDE \eqref{MV-2} in the sense of \autoref{martsol} with initial distribution $\mu_0$ on $\mD_T$.
\et

To prove the theorem above, we first prepare the following result about tightness.

\bl[Tightness]\label{lem:KM01}
The law of $(X_t^{N,1})_{t \in [0,T]}, N \in \mN$, in $\mD_T$ is tight.
\el

\begin{proof}
Let $A^N:=\int_0^\cdot(K*\eta_r^N)(X_r^{N,1}) \dif r$. By the boundedness of $K$, we have
$$
\sup_N \mE|A^N_t-A^N_s|^2\le\|K\|_\infty^2(t-s)^2,
$$
which, by the Kolmogorov-Chentsov criterion (cf. \cite[Theorem 23.7, p.511]{Ka3rd}), implies that $\{A^N\}$ is tight and the limiting processes are a.s continuous, and then $\{A^N\}$ is $C$-tight in $\mD_T$ (cf. \cite[Theorem 23.2 and Theorem 23.9]{Ka3rd}). Thus, the tightness of $X^{N,1}=A^N+L^{(\alpha),1}$ is a direct consequence of \cite[Corollary 3.33, p.\,353]{JS03}.
\end{proof}

Now we proceed to give the

\begin{proof}[Proof of \autoref{thm:wPC}]
Consider the following random measure with values in $\cP(\mD_T)$,
$$
\bar\omega \to \Pi_N(\bar\omega,\dif \omega) := \frac{1}{N}  \sum_{i=1}^N \delta_{X_\cdot^{N,i}(\bar\omega)}(\dif \omega) \in \cP (\mD_T),\quad \forall \bar\omega\in\Omega.
$$
By \autoref{lem:KM01} and \cite[ii) of Proposition 2.2, Chapter 1]{Szn91}, the laws of $\Pi_N$, $N \in \mN$, in $\cP( \cP (\mD_T))$ are tight. Without loss of generality, we assume that the laws of $\{\Pi_N\}_{N \in \mN}$ weakly converge to some $\Pi_\infty \in \cP(\cP(\mD_T))$. 

\medskip
Our aim below is to show that $\Pi_\infty$ is a Dirac measure, i.e.,
\begin{align}\label{CD01}
\Pi_\infty(\dif\zeta)=\delta_{\mu}(\dif\zeta),
\end{align}
where $\mu\in\sM^{b}_{\mu_0} $ is the unique martingale solution of DDSDE \eqref{MV-2} in the sense of \autoref{martsol} with initial distribution $\mu_0$. If we show the above assertion, then by \cite[i) of Proposition 2.2, Chapter 1]{Szn91}, we conclude \eqref{CC14}.

\medskip
To prove \eqref{CD01}, we adopt the classical martingale method (see \cite[Theorem 5.1]{HRZ24} or \cite{Zh23} for example). For given $f\in C^2_c(\mR^d)$ and $\zeta\in\cP(\mD_T)$, we define a functional on $\mD_T$ by
\begin{align}\label{eq:XM01-20241013}
M^{K}_{f,\zeta}(t,\omega):=f(w_t(\omega))-f(w_0(\omega))-\int^t_0 \sL_{K,\zeta_r}^{(\alpha)} f (w_r(\omega)) \dif r,\ \ \forall \omega \in \mD_T,
\end{align}
where $t\in[0,T]$, and $w_t: \mD_T \to \mR^d$ is the coordinate process, and 
$$
\sL_{K,\zeta_r}^{(\alpha)} f(x) := \sL^{(\alpha)} f (x) +((K*\zeta_r)\cdot\nabla f )(x) 
$$
with $\sL^{(\alpha)}$  defined by \eqref{eq:nonlocal}, and 
\begin{align}\label{eq:YU02}
\zeta_r:=\zeta\circ w_r^{-1}
\end{align}
is the marginal distribution of $\zeta$ at time $r\in[ 0,T]$. First of all, we claim that
\begin{itemize}
\item Claim 1: For $\Pi_\infty$-almost all $\zeta \in \cP(\mD_T)$, $\zeta_0= \mu_0$.
\item Claim 2: For $\Pi_\infty$-almost all $\zeta \in \cP(\mD_T)$, $(M^K_{f,\zeta})_{[0,T]}$ is a $(\cB_t)_{t \in [0,T]}$-martingale under $\zeta$.
\end{itemize}
If these two claims are established, then, by \autoref{martsol}, there is a $\Pi_\infty$-null set $\cN \subset \cP(\mD_T)$ such that  for all $\zeta\notin \cN$,
$$
\zeta\in\sM_{\mu_0}^b.
$$ 
Based on \autoref{thm:well}, $\sM_{\mu_0}^b$ only contains one point $\mu$, which derives that all the points $\zeta \notin \cN$ are equal to $\mu$ that is $\Pi_\infty (\{\mu\}) =1$ and then \eqref{CD01} is obtained. 

\medskip{
Now it remains to prove  Claim 1 and Claim 2. }

\medskip
\noindent
{\bf (Step 1)} For Claim 1, by the condition \eqref{CC15} and \cite[ii) of Proposition 2.2, Chapter 1]{Szn91}, one sees that
$$
\Pi_\infty\{\zeta\in\cP(\mD_T)\mid\zeta_0=\mu_0\}=1, 
$$
which is exactly Claim 1.

\medskip
\noindent
{\bf (Step 2)} To show Claim 2, we first give some notations and then divide the proof into two steps. Fix $n\in\mN$. Define
\begin{align*}
\mS^n_T:=\{ (s_1,s_2,...,s_n,s,t) \mid 0\le s_1\le s_2\le \cdot\cdot\cdot\le s_n\le s\le t\le T\}\subset [0,T]^{n+2}.
\end{align*}
For given  $g\in C_c(\mR^{nd})$, $h\in C([0,T]^{n+2})$ and $S=(s_1,s_2,...,s_n,s,t)\in \mS^n_T$, we also introduce functionals $\xi^g_f(S)$ and $\Xi^{g,h}_f$ on $\cP(\mD_T)$ by
\begin{align}\label{eq:XM02-20241013}
\xi^{g}_f(S,\zeta):=\int_{\mD_T} \(M^{K}_{f,\zeta}(t)-M^{K}_{f,\zeta}(s)\)(\omega) g(w_{s_1}(\omega),\cdots, w_{s_n}(\omega))\zeta(\dif \omega)
\end{align}
and
\begin{align}\label{eq:XM03-20241013}
\Xi^{g,h}_f(\zeta):=\int_{\mS^n_T}h(S)\xi^g_f(S,\zeta)\dif s_1\cdot\cdot\cdot\dif s_n\dif s\dif t.
\end{align}

\medskip
\noindent
{\bf (Step 2.1)}
In this step, we prove that for given $h\in C([0,T]^{n+2})$, $f\in C^2_c(\mR^d)$ and $g\in C_c(\mR^{nd})$,
\begin{align}\label{eq:YU01}
\int_{\cP(\mD_T)}|\Xi^{g,h}_f(\zeta)|\Pi_\infty(\dif \zeta)=\lim_{N\to\infty}\mE|\Xi^{g,h}_f(\Pi_N)|=0, \ \
\end{align}

\medskip
\noindent
$(1)$ We first show the second equality in \eqref{eq:YU01}. Notice that for any $\bar\omega\in\Omega$,
\begin{align*}
\big(\Pi_N(\bar\omega)\big)\circ w_s^{-1}(\dif x) =\frac1N\sum_{i=1}^N\delta_{X^{N,i}_s(\bar\omega)}(\dif x),
\end{align*}
which implies
\begin{align*}
\xi^{g}_f(S,\Pi_N) = & \frac{1}{N}\sum_{i=1}^N\(M_{f,\Pi_N}^K(t,X^{N,i}_\cdot) -M_{f,\Pi_N}^K(s,X^{N,i}_\cdot)\)
g\(X^{N,i}_{s_1},\cdots, X^{N,i}_{s_n}\), \ \  \mP-a.s.,
\end{align*}
and implies, by applying It\^o's formula (cf. \cite{IW89}) to $f(X_t^{N,i})$ with $f\in C^2_c(\mR^d)$,  that for all most surely $\bar \omega \in \Omega$,
\begin{align*}
M_{f,\Pi_N(\bar\omega)}^K(t,X^{N,i}_\cdot (\bar\omega))= \left( \int_0^t \int_{\mR^d}\left(f(X^{N,i}_{r-}+z)-f(X^{N,i}_{r-})\right)\widetilde{N}^i(\dif z,\dif r)\right)(\bar \omega).
\end{align*}
Thus, the desired result follows from the definition \eqref{eq:XM03-20241013}, and the fact,
\begin{align*}
\mE|\xi^{g}_f(S,\Pi_N)|^2&=\frac{1}{N^2}\mE\left|\sum_{i=1}^N\int_s^t \int_{\mR^d}\left(f(X^{N,i}_{r-}+z)-f(X^{N,i}_{r-})\right)g\(X^{N,i}_{s_1},\cdots, X^{N,i}_{s_n}\)\widetilde{N}^i(\dif z,\dif r)\right|^2\\
&=\frac{1}{N^2}\sum_{i=1}^N\mE\int_s^t \int_{\mR^d}\left|\left(f(X^{N,i}_{r-}+z)-f(X^{N,i}_{r-})\right)g\(X^{N,i}_{s_1},\cdots, X^{N,i}_{s_n}\)\right|^2 \nu^{(\alpha)}(\dif z)\dif r\\
&\overset{\eqref{eq:XM02}}{\lesssim} N^{-1}\|f\|_{\bC^1}^2\|g\|_\infty^2,
\end{align*}
where we used the independence, and It\^o's isometry (cf. \cite[Section 3 of Chapter II]{IW89}), and the fact $$g\(X^{N,i}_{s_1},\cdots, X^{N,i}_{s_n}\)\in \sF_{s}$$ in the second equality.
 
\medskip\noindent
$(2)$ As for the first equality in \eqref{eq:YU01}, it follows from the fact 
\begin{align}\label{CD00}
\Xi^{g,h}_f\in C_b(\cP(\mD_T)),
\end{align}
and the weak convergence of $\Pi_\infty$ and the laws of $\{\Pi_N\}_{N \in \mN}$. Now it remains to prove \eqref{CD00}. Since $K \in C_b(\mR^d)$ and the boundness of the operator $\sL^{(\alpha)}$ (see \cite[Lemma 4.1]{HWW23}), one sees that the functional $\Xi_f^{g,h}$ is bounded by the definitions \eqref{eq:XM01-20241013}, \eqref{eq:XM02-20241013}, and \eqref{eq:XM03-20241013}. Notice that   $\Xi_f^{g,h}$ is a nonlinear functional of $\zeta$, we have to take some care for its continuity.

Suppose that $\zeta^{(m)} \in \cP (\mD_T )$ weakly converges to $\zeta  \in \cP (\mD_T )$. Based on Skorokhod's representation theorem, there is a probability space 
$$
(\widetilde{\Omega},\widetilde{\sF},(\widetilde{\sF}_s)_{s\ge0},\widetilde{\mP})
$$
and a sequence of c\'adl\'ag processes $\{\widetilde{X}^{(m)}\}_{m=1}^\infty$ as well as another c\'adl\'ag process $\widetilde{X}$ thereon such that 
\begin{align}\label{eq:XM05-20241013}
    \widetilde{\mP}\circ(\widetilde{X}^{(m)})^{-1}=\zeta^{(m)},\quad \widetilde{\mP}\circ(\widetilde{X})^{-1}=\zeta,
\end{align}
and
\begin{align*}
    \lim_{m\to\infty}\widetilde{X}^{(m)}=\widetilde{X} \quad \text{in the Skorokhod topology}, \quad \widetilde{\mP}-a.s.
\end{align*}
We claim that
\begin{align}\label{eq:XM04-20241023}
|\xi_f^{g}(S,\zeta^{(m)})  - \xi_f^{g} & (S,\zeta)| \lesssim \sum_{u\in S}\widetilde{\mE}\left[|\widetilde{X}^{(m)}_u-\widetilde{X}_u|^\beta\wedge 1\right]+\int_s^t \widetilde{\mE}\left[|\widetilde{X}^{(m)}_r-\widetilde{X}_r|^\beta\wedge1\right] \dif r.
\end{align}
If we get this claim, then by the definition \eqref{eq:XM03-20241013} and $\dif S:=\dif s_1\cdot\cdot\cdot\dif s_n\dif s\dif t$, we have
\begin{align*}
    |\Xi_f^{g,h}(\zeta^{(m)})-\Xi_f^{g,h}(\zeta)|&\le \|h\|_\infty \int_{\mS^n_T}|\xi_f^{g}(S,\zeta^{(m)})-\xi_f^{g}(S,\zeta)|\dif S\\
    &\lesssim \int_{\mS^n_T}\left(\sum_{u\in S}\widetilde{\mE}\left[|\widetilde{X}^{(m)}_u-\widetilde{X}_u|^\beta\wedge 1\right]+\int_s^t \widetilde{\mE}\left[|\widetilde{X}^{(m)}_r-\widetilde{X}_r|^\beta\wedge1\right] \dif r\right)\dif S.
\end{align*}
Furthermore, observing that, by \cite[2.3 of Chapter VI, p.339]{JS03}, for $\widetilde{\mP}$ almost surely $\tilde{\omega}\in\widetilde{\Omega}$, $\lim_{m\to\infty}\widetilde{X}^{(m)}_r(\tilde{\omega})=\widetilde{X}_r(\tilde{\omega})$ for all 
$$
r\in \{u\in[0,T]\mid \widetilde{X}_{u-}(\tilde{\omega})=\widetilde{X}_{u}(\tilde{\omega})\},
$$
which is dense in $[0,T]$ (cf. \cite[1.7 of Chapter VI, p.326]{JS03}), and using the dominated convergence theorem, we have 
\begin{align*}
    \lim_{m\to\infty}|\Xi_f^{g,h}(\zeta^{(m)})-\Xi_f^{g,h}(\zeta)|=0
\end{align*}
and then obtain the claim \eqref{eq:YU01}.

Now we proceed to prove the claim \eqref{eq:XM04-20241023}. By the definitions, we have 
\begin{align*}
|\xi_f^{g}(S,\zeta^{(m)})  - \xi_f^{g} & (S,\zeta)| \leq 
 \left |\int_{\mD_T} \[f(w_t(\omega)) - f(w_s(\omega)) \] g(w_{s_1}(\omega), \cdots, w_{s_n}(\omega)) (\zeta^{(m)} - \zeta) (\dif \omega)\right|\\
& + \left |\int_{\mD_T} \[ \int_s^t \sL_{K,\zeta_r}^{(\alpha)}f (w_r(\omega)) \dif r\] g(w_{s_1}(\omega), \cdots, w_{s_n}(\omega)) (\zeta^{(m)} - \zeta) (\dif \omega)\right|\\
&+\| \nabla f \|_\infty\|g\|_\infty \int_s^t \int_{\mD_T}  |K* (\zeta^{(m)} - \zeta)_r |(w_r(\omega)) \zeta^{(m)}(\dif \omega)  \dif r\\
=:& \sJ_m^{(1)} + \sJ_m^{(2)}  + \| \nabla f \|_\infty\|g\|_\infty\sJ_m^{(3)} .
\end{align*}
By \eqref{eq:XM05-20241013} and \eqref{eq:YU02}, we have
\begin{align*}
 \sJ_m^{(1)}\le \|f\|_{\bC^1_b}\|g\|_{\bC^1_b}\sum_{u\in S}\widetilde{\mE}\left[|\widetilde{X}^{(m)}_u-\widetilde{X}_u|\wedge 1\right],
\end{align*}
and,
\begin{align*}
\sJ_m^{(2)}
& \lesssim (\|\sL^{(\alpha)}f\|_{\bC^\beta}+\|K\|_{\bC^\beta}\|\nabla f\|_{\bC^\beta})\|g\|_{\bC^\beta}\sum_{u\in S}\widetilde{\mE}\left[|\widetilde{X}^{(m)}_u-\widetilde{X}_u|^\beta\wedge 1\right]\\
& \lesssim \sum_{u\in S}\widetilde{\mE}\left[|\widetilde{X}^{(m)}_u-\widetilde{X}_u|^\beta\wedge 1\right],
\end{align*}
where we used \cite[Lemma 4.1]{HWW23} again in the last inequality. Moreover, by \eqref{eq:XM05-20241013},
\begin{align*}
\sJ_m^{(3)}&=\int_s^t \int_{\mD_T}  |\widetilde{\mE}K(w_r(\omega)-\widetilde{X}^{(m)}_r)-\widetilde{\mE}K(w_r(\omega)-\widetilde{X}_r)| \zeta^{(m)}(\dif \omega)  \dif r\\
    &\le \|K\|_{\bC^\beta} \int_s^t \int_{\mD_T} \widetilde{\mE}\left[|\widetilde{X}^{(m)}_r-\widetilde{X}_r|^\beta\wedge1\right] \zeta^{(m)}(\dif \omega)  \dif r\\
& \lesssim \int_s^t \widetilde{\mE}\left[|\widetilde{X}^{(m)}_r-\widetilde{X}_r|^\beta\wedge1\right] \dif r.
\end{align*}
Combining the calculations above, we get the claim \eqref{eq:XM04-20241023}.

\medskip
\noindent
{\bf (Step 2.2)}
Then it follows from \eqref{eq:YU01} that for any function $h\in C([0,T]^{n+2})$, $f\in C^2_c(\mR^d)$ and $g\in C_c(\mR^{nd})$
$$
\int_{\cP(\mD_T)}|\Xi_f^{g,h}(\zeta)| \Pi_\infty(\dif \zeta) =0.
$$
Hence, observing that $C([0,T]^{n+2})$, $C^2_c(\mR^d)$ and $C_c(\mR^{nd})$ are separable, one can find a common $\Pi_\infty$-null set $\cN_1\subset\cP(\mD_T)$
such that for all $\zeta\notin\cN_1$ and for all 
$f\in C_c^2(\mR^{d})$, $g\in C_c(\mR^{nd})$, $h\in C([0,T]^{n+2})$, 
\begin{align*}
 \Xi^{g,h}_f(\zeta)\overset{\eqref{eq:XM03-20241013}}{=}\int_{\mS^n_T}\xi^{g}_f(S,\zeta)h(S)\dif S=0,
\end{align*}
which implies that for Lebegue almost surely $S=(s_1,s_2,..,s_n,s,t)\in \mS^n_T$,
$$
\xi^{g,h}_f(S,\zeta)\overset{\eqref{eq:XM02-20241013}}{=}\int_{\mD_T} \(M^{K}_{f,\zeta}(t)-M^{K}_{f,\zeta}(s)\)(\omega) g(w_{s_1}(\omega),\cdots, w_{s_n}(\omega))\zeta(\dif \omega)=0.
$$Recalling the definition \eqref{eq:XM01-20241013}, since $t\to w_t(\omega)$ is right continuous, we further have that for every $\zeta\notin\cN_1$ and for any $S\in \mS^n_T$,
$$
\int_{\mD_T} \(M^{K}_{f,\zeta}(t)-M^{K}_{f,\zeta}(s)\)(\omega) g(w_{s_1}(\omega),\cdots, w_{s_n}(\omega))\zeta(\dif \omega)=0,
$$
which implies that $M^K_{f,\zeta}$ is a $(\cB_t)_{t \in [0,T]}$-martingale under $\zeta$ and then the proof of Claim 2 is complete.
\end{proof}

\subsubsection{Strong propagation of chaos}\label{subsubsec:S-chaos}

In this subsubsection, under $\beta> 1-\alpha/2$, we use Zvonkin's transform and the previous weak convergence result \autoref{thm:wPC} to derive the strong propagation of chaos for the interacting $N$-particle system \eqref{eq:AC01}. To this end, we consider the following coupled DDSDE:
\begin{align}\label{eq:XM04}
\dif X_t^{j} = (K*\mu^j_t)(X_t^{j}) \dif t + \dif L^{(\alpha),j}_t,\ \ X_0^{j}=\xi_j,
\end{align}
where $\mu^j_t$ is the time marginal distribution of the solution $X_t^{j}$ and $\{\xi_j\}_{j=1}^\infty$ is a family of i.i.d. random variables. By the uniqueness in \autoref{thm:well}, one sees that $\mu_t^j, j \in \mN$ are all the same.
 
\bt[Strong propagation of chaos]\label{thm:SPC}
Let $\alpha \in (0,1)$. Assume that $K\in\bC^\beta(\mR^d)$ with some $\beta \in(1-\alpha/2,1)$, and $(\xi_1^N,\cdots,\xi^N_N,\cdots)=(\xi_1,\cdots,\xi_N,\cdots)$. Then for every $T>0$ and $p>0$,
\begin{align}\label{S:CC14}
\lim_{N\to\infty}\sup_{j=1,..,N}\mE\(\sup_{t\in[0,T]}|X^{N,j}_t-X^j_t|^p\)=0.
\end{align}
\et

\begin{proof}
Similar with the proof of \autoref{thm:wEuler}, we use Zvonkin's transform to prove this theorem. It suffices to prove the case $p\geq 2$ since Jensen's inequality. First of all, define
\begin{align}\label{eq:DT01}
\tilde b_{\eta^N}(t,x):= (K*\eta_t^N)(x),\ \ \text{and}\ \ \tilde b_\mu(t,x):= (K*\mu_t)(x),
\end{align}
where $\mu_t = \mu_t^j$, $j=1,\cdots,N$. Let $u$ be the unique solution of PDE \eqref{eq:XM00} with the drift coefficient $\tilde b_\mu(t,x):= (K*\mu_t)(x)$ and the large enough positive number $\lambda$, and $\widetilde{N}^j(\dif s,\dif z)$ be the compensated Poisson random measure to $\alpha$-stable process $L^{(\alpha),j}$. Denoting by $\Phi (t,x):= x + u(t,x)$, and
$$
M_t^j := \int_0^t \int_{\mR^d}
\[ u (s,X^j_{s-}+z) -  u (s,X^j_{s-}) \]\widetilde{N}^j(\dif s,\dif z) + \int_0^t \int_{\mR^d} z {N}^j(\dif s,\dif z),
$$
and
$$
M_t^{N,j} := \int_0^t \int_{\mR^d}
\[ u(s,X^{N,j}_{s-}+z) - u(s,X^{N,j}_{s-}) \]\widetilde{N}^j(\dif s,\dif z)+ \int_0^t \int_{\mR^d} z {N}^j(\dif s,\dif z),
$$
and 
\begin{align*}
R^{N,j}_t:=\int_0^t \left( ( \tilde b_{\eta^N}-  \tilde b_\mu )\cdot ( \mathbb{I} + \nabla  u) \right)(s,X^{N,j}_s)\dif s,
\end{align*}
and applying It\^o's formula with SDEs \eqref{eq:AC01} and \eqref{eq:XM04}, one sees that,  by PDE \eqref{eq:XM00}, 
\begin{align*}
\Phi(t,X^{j}_t)=\Phi(0,\xi_j)+\lambda\int_0^t u(s,X^{j}_s)\dif s+M_t^j
\end{align*}
and
\begin{align*}
\Phi(t,X^{N,j}_t)&=\Phi(0,\xi^N_j)+\lambda\int_0^t u(s,X^{N,j}_s)\dif s+R^{N,j}_t +M_t^{N,j}.
\end{align*}
Hence, by the definitions and H\"older's inequality, we have that for $p\geq 2$,
\begin{align*}
\mE\left[\sup_{s\in[0,t]}|X^j_s-X^{N,j}_s|^p\right]
& \lesssim_p \|\nabla u\|_{\mL^\infty_T}^p \mE\left[\sup_{s\in[0,t]}|X^j_s-X^{N,j}_s|^p\right]+ \mE\left[\sup_{s\in[0,t]} |M_s^j - M^{N,j}_s|^p \right] \\
+ &  \lambda^p T^{p-1}\|\nabla u\|_{\mL^\infty_T}^p\int_0^t  \mE\left[\sup_{r\in[0,s]}|X^j_r-X^{N,j}_r|^p\right] \dif s
  + \mE\left[\sup_{s\in[0,t]} \mE  |R^{N,j}_s|^p\right], 
\end{align*}
which, by \eqref{eq:XM01} and the same tricks as the proof of  \eqref{eq:XM05}, we get that
\begin{align}\label{eq:XM06}
\mE\left[\sup_{s\in[0,t]}|X^j_s-X^{N,j}_s|^p\right] \lesssim_{\lambda,p,T}   \int_0^t  \mE\left[\sup_{r\in[0,s]}|X^j_r-X^{N,j}_r|^p\right] \dif s + \mE\left[\sup_{s\in[0,t]} \mE  |R^{N,j}_s|^p\right].
\end{align}
Now we claim that for any $p\geq 2$,
\begin{align}\label{S5:10}
\lim_{N\to\infty}\sup_{1\leq j \leq N}\mE\left[\sup_{s\in[0,t]}|R^{N,j}_s|^p\right]=0,
\end{align}
which, combining \eqref{eq:XM06} with Gronwall's inequality \autoref{lem:g}, conculdes \eqref{S:CC14}. 

\medskip
Thus, it remains to show \eqref{S5:10}. Indeed, one sees that by \eqref{eq:DT01} and \eqref{eq:DT02},
\begin{align}\label{eq:XM09}
\mE\left[\sup_{s\in[0,t]}|R^{N,j}_s|^p\right]& \lesssim_{T,p} \frac{\|\nabla\Phi\|^p_{\mL^\infty_T}}{N^p} \int_0^t \mE\left|\sum_{k=1}^N\left(K(X^{N,j}_s-X^{N,k}_s)-\tilde b_{\mu_s} (X^{N,j}_s)\right)\right|^p\dif s \nonumber\\
& \lesssim \frac{ \|K\|_\infty^{p-2}}{N^{2}} \int_0^t \mE\left|\sum_{k=1}^N\left(K(X^{N,j}_s-X^{N,k}_s)-\tilde b_{\mu_s} (X^{N,j}_s)\right)\right|^{2}\dif s.
\end{align}
Without loss of generality, we assume that $N >2$. Define
$$
Z_s^{N,j,k}:=K(X^{N,j}_s-X^{N,k}_s)-\tilde b_{\mu_s} (X^{N,j}_s)
$$
and the set
\begin{align*}
J := \Big\{(k_1, &\, k_2)\in \mN^2 \mid 1\leq k_1,k_2\leq N, \text{ and } k_1\neq k_2\neq j \Big\}.
\end{align*}
Notice that, since $X_s^{1}, X_s^{2}, X_s^{3}$  are independent and have the same distribution $\mu_s$, by defnitions and \eqref{CC14} and \eqref{eq:DT01}, we have
\begin{align}\label{eq:XM08}
\lim_{N \to \infty} \mE\< Z_s^{N,3,1},Z_s^{N,3,2} \>   =0 .
\end{align}
Moreover, one sees that
\begin{align}\label{eq:XM10}
&\mE\left|\sum_{k=1}^N\left(K(X^{N,j}_s -X^{N,k}_s)   -\tilde b_{\mu_s} (X^{N,j}_s)\right)\right|^{2}= \mE\left|\sum_{k=1}^N Z_s^{N,j,k}\right|^{2}\nonumber \\
&\qquad = \mE\left [\left(\sum_{J^{c}}+ \sum_{J}\right) \< Z_s^{N,j,k_1},Z_s^{N,j,k_2}\>\right]\nonumber \\
& \qquad
 \leq 4 N  \|K\|_\infty^2 + (N-1)(N-2) \mE\< Z_s^{N,3,1},Z_s^{N,3,2}\>,
\end{align}
where, in the last inequality, we used the fact that for every element in $J$,
$$ (X^{N,k_1},X^{N,k_2},X^{N,j})\overset{(d)}{=}(X^{N,1},X^{N,2},X^{N,3}).
$$
Hence, back to \eqref{eq:XM09}, combining with \eqref{eq:XM08} and \eqref{eq:XM10}, one sees that
\begin{align*}
\lim_{N \to \infty }\sup_{1\leq j \leq N}\mE\left[\sup_{s\in[0,t]}|R^{N,j}_s|^p\right] \lesssim \lim_{N \to \infty } \frac{1}{N}+ 0 =0.
\end{align*}
Thus we obtain \eqref{S5:10} and the proof is completed. 
\end{proof}

\subsection{Propagation of chaos for Euler's scheme}
Below, fix $h\in(0,1)$, and let $\{\xi_i\}_{i=1}^\infty$ be a sequence of i.i.d. random variables
in $\mR^d$ with the common distribution $\mu_0 \overset{(d)}{=} X_0$, $X^{N,i,h}$ be the unique solution of the following Euler's scheme for the $N$-particle systems \eqref{eq:AC01} for each $1\leq i\leq N$ with fixed $N\in\mN$: \begin{align}\label{eq:YT00}
\dif X^{N,i,h}_t=\(K*\eta^{N,h}_{\pi_h(t)}\)(X^{N,i,h}_{\pi_h(t)})\dif t+\dif L^{(\alpha),i}_t,\quad X^{N,i,h}_0\overset{(d)}{=} \xi_i,
\end{align}
where $\pi_h(t):=[t/h]h$, and $\{L^{(\alpha),i}\}_{i\geq 1}$ be a sequence of independent $d$-dimensional $\alpha$-stable processes with $\alpha \in (0,1)$, and $\eta^{N,h}_t$ is the empirical measure of $\{X^{N,i,h}_t, ~i=1,..,N\}$ defined by \begin{align*}
\eta^{N,h}_t(\dif y):= \frac{1}{N} \sum_{j=1}^N \delta_{X_t^{N,j,h}}(\dif y),
\end{align*}
where $\delta_x$ stands for the Dirac measure concentrated at point $x\in \mR^d$. In the following, for each $i\in \mN$, let $X^{i,h}$ be the unique  solution of the following coupled Euler's scheme:
\begin{align}\label{eq:YT01}
\dif X^{i,h}_t=\(K*\mu^{i,h}_{\pi_h(t)}\)(X^{i,h}_{\pi_h(t)})\dif t+\dif \widehat L^{(\alpha),i}_t,\quad  X^{i,h}_0 \overset{(d)}{=} \xi_i,
\end{align}
where $\mu^{i,h}_t$ is the distribution of $X^{i,h}_t$, and $\widehat L^{(\alpha),i}$ is also a sequence of i.d.d. $\alpha$-stable processes. Notice that $\{X^{i,h}\}_{i=1}^\infty$ is a family of i.i.d. stochastic processes with common distribution as $X^{h}$.

Here is the main result in this subsection. 

\bt[Propagation of chaos for Euler's scheme]\label{thm:PCE}
Let $\alpha \in (0,1)$. Assume $K\in \bC^\beta(\mR^d)$ with some $\beta\in(0,1)$. Then for any $h\in(0,1)$ and $T>0$, we have that 
\begin{enumerate}[\rm (i)]
\item if $X^{i,h}_0= X^{N,i,h}_0$ and $L^{(\alpha),i} = \widehat L^{(\alpha),i}$ for each $i$, then for every $p>0$,
\begin{align}\label{S7:00}
\lim_{N\to\infty}\sup_{i=1,..,N}\mE\left(\sup_{t\in[0,T]}|X^{N,i,h}_t-X^{i,h}_t|^p\right)=0;
\end{align}
\item letting $\mu^{N,1,h}_t$ be the time marginal distribution of $(X^{N,1,h}_t)_{t\in[0,T]}$, 
\begin{align*}
\lim_{N\to\infty}\sup_{t\in[0,T]}\|\mu^{N,1,h}_t-\mu^{1,h}_t\|_{\beta,\var}=0.
\end{align*}
\end{enumerate}
\et
\begin{proof}

(i) The proof is based on the technique used in \cite[Theorem 1.3]{Zh19}. By Jensen's inequality, it suffices to show the case $p\ge 2$. First of all, we claim that there is a positive constant $c$ independent of $N$, such that
\begin{align}\label{S6:01}
\sup_{1\leq i\leq N}\mE\left(\sup_{t\in[0,T]}|X^{N,i,h}_t-X^{i,h}_t|^p\right)\lesssim_c  \int_0^T \left (\sup_{1\leq i\leq N} \mE|X^{N,i,h}_{\pi_h(s)}-X^{i,h}_{\pi_h(s)}|^{2\beta}\right)\dif s+1/N.
\end{align}
From this, by the induction method, it is easy to derive the desired result \eqref{S7:00}. Indeed, we have that, by the assumption,
\begin{align*}
\sup_{1\leq i \leq N}\mE|X^{N,i,h}_0-X^{i,h}_0|^p=0,
\end{align*}
and if we suppose that the following holds for some $n\geq 1$,
\begin{align*}
\lim_{N\to\infty}\sup_{1\leq i \leq N}\mE\left(\sup_{t\in[0,nh]}|X^{N,i,h}_t-X^{i,h}_t|^p \right)=0,
\end{align*}
then, observing the fact,
\begin{align*}
\sup_{1\leq i \leq N}\mE\left(\sup_{t\in[0,(n+1)h]}|X^{N,i,h}_t-X^{i,h}_t|^p\right) \overset{\eqref{S6:01}}{\lesssim} \sum_{k=0}^n \sup_{1\leq i \leq N}\mE|X^{N,i,h}_{kh}-X^{i,h}_{kh}|^{2\beta} +1/N,
\end{align*}
and the induction
hypothesis with $\beta \in (0,1)$, we obtain that
\begin{align*}
\lim_{N\to\infty}\sup_{1\leq i \leq N}\mE\left(\sup_{t\in[0,(n+1)h]}|X^{N,i,h}_t-X^{i,h}_t|^p \right)=0.
\end{align*}

\medskip
Now we show the inequality \eqref{S6:01}. By definitions, and SDEs \eqref{eq:YT00} and \eqref{eq:YT01}, we have that for $p\ge 2$, 
\begin{align*}
&\, \mE\left(\sup_{t\in[0,T]}|X^{N,i,h}_t   -X^{i,h}_t|^p\right) 
\lesssim \|K\|_\infty^{p-2} \mE\left(\sup_{t\in[0,T]}|X^{N,i,h}_t-X^{i,h}_t|^2\right)\\
& \qquad \lesssim  \mE \int_0^T \left|\(K*\mu^{h}_{\pi_h(s)}\)(X^{i,h}_{\pi_h(s)})-\(K*\widetilde\eta^{N,h}_{\pi_h(s)}\)(X^{i,h}_{\pi_h(s)})\right|^2\dif s\\
&\qquad +\mE \int_0^T \left|\(K*\eta^{N,h}_{\pi_h(s)}\)(X^{N,i,h}_{\pi_h(s)})-\(K*\widetilde\eta^{N,h}_{\pi_h(s)}\)(X^{i,h}_{\pi_h(s)})\right|^2 \dif s\\
& \qquad =: \sI_1^{N,i,h}+\sI_2^{N,i,h},
\end{align*}
where 
\begin{align*}
\widetilde\eta^{N,h}_t(\dif y):=\frac1N\sum_{j=1}^N\delta_{X^{j,h}_t}(\dif y).
\end{align*}
Let us first treat the term $\sI_1^{N,i,h}$. Denoting 
$$
Z_{s}^{h,i,j}:=\(K*\mu^{h}_{\pi_h(s)}\)(X^{i,h}_{\pi_h(s)})-K(X^{i,h}_{\pi_h(s)}-X^{j,h}_{\pi_h(s)}),
$$
by the independence of $\{X^{i,h}_\cdot\}_{i=1}^N$, we have that for any $i\ne j\ne k$, 
\begin{align}\label{eq:YT03}
\mE \<Z_s^{h,i,j}, Z_s^{h,i,k}\> =0,
\end{align}
where we used the fact that $\mu^h_{\pi_h(s)}$ is the law of $X^{j,h}_{\pi_h(s)}$ for any $j\in \mN$. Consequently, we obtain that
\begin{align*}
&\mE\left|\(K*\mu^{h}_{\pi_h(s)}\)(X^{i,h}_{\pi_h(s)})-\(K*\widetilde \eta^{N,h}_{\pi_h(s)}\)(X^{i,h}_{\pi_h(s)})\right|^2\\
=&\frac{1}{N^2}\mE\left|\sum_{j=1}^N Z_s^{h,i,j}\right|^2
= \frac1{N^2}\sum_{j,k=1}^N \mE \< Z_s^{h,i,j}, Z_s^{h,i,k}\>
\\
\overset{\eqref{eq:YT03}}{\leq} & \frac1{N^2}\left(2\sum_{k=1}^N \mE |\<Z_s^{h,i,i},Z_s^{h,i,k}\>|+\sum_{j=1}^N \mE |Z_s^{h,i,j}|^2\right)\lesssim \frac{\|K\|_\infty^2}{N},
\end{align*}
which implies that
\begin{align}\label{S6:00}
\sup_{1\leq i \leq N}\sI_1^{N,i,h}\lesssim 1/N.
\end{align}
For $\sI_2^{N,i,h}$, we have
\begin{align*}
\sI_2^{N,i,h} &=\frac1{N^2}\mE \int_0^T \left|\sum_{j=1}^N\left(K(X^{N,i,h}_{\pi_h(s)}-X^{N,j,h}_{\pi_h(s)})-K(X^{i,h}_{\pi_h(s)}-X^{j,h}_{\pi_h(s)})\right)\right|^2\dif s\\
&\lesssim \frac1{N}\sum_{j=1}^N\mE \int_0^T \left|\left(K(X^{N,i,h}_{\pi_h(s)}-X^{N,j,h}_{\pi_h(s)})-K(X^{i,h}_{\pi_h(s)}-X^{j,h}_{\pi_h(s)})\right)\right|^2 \dif s\\
&\lesssim \|K\|_{\bC^{\beta}}^2  \int_0^T \left(\sup_{1\leq j \leq N}\mE |X^{N,j,h}_{\pi_h(s)}-X^{j,h}_{\pi_h(s)})|^{2\beta}\right)\dif s,
\end{align*}     
which, combining with \eqref{S6:00}, implies \eqref{S6:01}. The proof is completed.

\medskip\noindent
(ii) Without loss of generality, we assume that the sequence, consisting of $X^{N,i,h}_0$,  $  X_0^{i,h}$, $L^{(\alpha),i}$, $\widehat L^{(\alpha),i}$, $i \in \mN$, is independent (of course, defined on the same probability space as well). Letting
$\bL^{(\alpha),N} := (L^{(\alpha),1},\cdots, L^{(\alpha),N})$, and $\widehat \bL^{(\alpha),N} := (\widehat  L^{(\alpha),1},\cdots, \widehat  L^{(\alpha),N})$, and 
$$
\bX^{N,N,h} := (X^{N,1,h},\cdots, X^{N,N,h}),\quad\text{and}\quad  \bX^{N,h} := (X^{1,h},\cdots, X^{N,h}),
$$
one sees that the Euler's scheme \eqref{eq:YT00} and \eqref{eq:YT01}, respectively, also be solved by SDEs in $\mR^{dN}$: 
\begin{align}\label{eq:DTT-1}
\dif \bX_t^{N,N,h} = \cK^N(\bX_{\pi_h(t)}^{N,N,h})\dif t + \dif \bL_t^{(\alpha),N},\qquad \bX_0^{N,N,h} \overset{(d)}{=} (\xi_1,\cdots, \xi_N),
\end{align}
and
\begin{align}\label{eq:DTT-2}
\dif \bX_t^{N,h} = ({ \boldsymbol{K}^N*\boldsymbol{\mu}_{\pi_h(t)}^{N,h})( \bX_{\pi_h(t)}^{N, h})}\dif t + \dif \widehat\bL_t^{(\alpha),N},\qquad \bX_0^{N,h} \overset{(d)}{=} (\xi_1,\cdots, \xi_N),
\end{align}
where 
and $\cK^N$ is defined by \eqref{eq:JU93}, and 
$
\boldsymbol{K}^N*\boldsymbol{\mu}_{t}^{N,h}:= ({K}*{\mu}_{t}^{1,h}, \cdots, {K}*{\mu}_{t}^{N,h}).
$

By the definition of Euler's scheme \eqref{eq:DTT-2} and \eqref{eq:DTT-1}, it is easy to see that there are two Borel functions $f_{N,N},g_N:\mR^{Nd}\times \mD_T^{\times N} \to \mD_T^{\times N}$ such that 
$$
\bX^{N,N,h} =f_{N,N}(\bX_0^{N,N,h}, \bL^{(\alpha),N}),\qquad 
 \bX^{N,h}=g_N( \bX_0^{N,h}, \widehat \bL^{(\alpha),N}).
$$
Considering the process $ \bY^{N,N,h}=(Y^{N,1,h},...,Y^{N,N,h})$ defined by
$$
\bY^{N,N,h}:=g_N(\bX_0^{N,N,h},  \bL^{(\alpha),N}),
$$
one sees that, by $(\bX_0^{N,N,h} , \bL^{(\alpha),N} ) \overset{(d)}{=} ( \bX_0^{N,h}, \widehat \bL^{(\alpha),N})$, we have
\begin{align}\label{eq:XM00-20241015}
\bY^{N,N,h} \overset{(d)}{=} \bX^{N,h},
\end{align} 
and then, by definitions, 
\begin{align}
\sup_{t\in[0,T]}\|\mu^{N,1,h}_t-\mu^{1,h}_t\|_{\beta,\var}&=\sup_{t\in[0,T]}\sup_{\|\varphi\|_{\bC^\beta}\le 1}\left|\mE\varphi(  X^{N,1,h}_t)-\mE\varphi(  X^{1,h}_t)\right|\nonumber\\
&\overset{\eqref{eq:XM00-20241015}}{=}\sup_{t\in[0,T]}\sup_{\|\varphi\|_{\bC^\beta}\le 1}\left|\mE\( \varphi(  X^{N,1,h}_t)-\varphi( Y^{N,1,h}_t) \)\right|\nonumber\\
&\le \sup_{t\in[0,T]}\mE\left| X^{N,1,h}_t- Y^{N,1,h}_t\right|^\beta. \label{eq:XM00-20241016}
\end{align}
Moreover, by definitions, it is easy to see that  $\bY^{N,N,h}_0=\bX^{N,N,h}_0$ and $\bY^{N,N,h}$ solves \eqref{eq:DTT-2} driven by the process $\bL^{{\alpha},N}$. Hence, by the proof of \eqref{S7:00}, we obtain that
$$
\lim_{N \to \infty} \sup_{t\in[0,T]}\mE\left| X^{N,1,h}_t- Y^{N,1,h}_t\right|^\beta =0,
$$
which together with \eqref{eq:XM00-20241016} derives the desired result. The proof is finished.  
\end{proof}

\medskip

\subsection*{Acknowledgments}
We are deeply grateful to Prof. Xicheng Zhang for his valuable suggestions and for correcting some errors. 

\end{document}